\documentclass[
	framed_theorems,
	framed_proofs,
	print,
]{OCPreprint}

\usepackage[ruled,vlined,linesnumbered,algo2e,algosection]{algorithm2e}

\usepackage{etoolbox}
\makeatletter
\patchcmd{\@algocf@start}{\addtolength{\hsize}{-1.5em}}{}{}{}
\makeatother

\usepackage{pgfplots}
\usepgfplotslibrary{fillbetween}
\usepackage{pgfplotstable}
\pgfplotsset{compat=1.16}

\usepackage{xcolor}

\usepackage{enumitem}
\setlist[enumerate,1]{label=\textup{(\roman*)}}

\AtEndPreamble{
	\newtheorem{notation}[theorem]{Notation}
	\ifthenelse{\boolean{framedtheorems}}{
		\maketheoremcolored{notation}
		\colorlet{notationframecolor}{assumptionframecolor}
		\colorlet{notationshadecolor}{assumptionshadecolor}
	}{}
}

\newcommand\TV{\operatorname{TV}}
\newcommand\BV{\operatorname{BV}}
\newcommand\F{\mathbb{F}}

\newcommand\STfull{\hyperref[eq:ST]{\textup{(ST($\hat n, \hat a$))}}}
\definecolor{ucol}{RGB}{255,128,0}
\definecolor{Kucol}{RGB}{0,0,250}
\definecolor{gruen}{RGB}{0,153,0}
\definecolor{fcol}{RGB}{255,40,20}
\definecolor{nablacol}{RGB}{51,153,255}

\title[Integer optimal control problems with total variation regularization]%
{Integer optimal control problems with total variation regularization: Optimality conditions and fast solution of subproblems}
\author[Marko, Wachsmuth]{%
	Jonas Marko%
	\footnote{%
	Brandenburgische Technische Universität Cottbus--Senftenberg,
	Institute of Mathematics,
	03046 Cottbus,
	Germany,
	\email{markojo1@b-tu.de},
	\url{https://www.b-tu.de/fg-optimale-steuerung}
	},
	Gerd Wachsmuth%
	\footnote{%
	Brandenburgische Technische Universität Cottbus--Senftenberg,
	Institute of Mathematics,
	03046 Cottbus,
	Germany,
	\email{wachsmuth@b-tu.de},
	\url{https://www.b-tu.de/fg-optimale-steuerung},
	ORCID: \href{https://orcid.org/0000-0002-3098-1503}{0000-0002-3098-1503}%
	}
}

\begin{document}

\maketitle

 \begin{abstract}
	We investigate local optimality conditions of first and second order for integer optimal control problems with total variation regularization via a finite-dimensional switching point problem.
	We show the equivalence of local optimality for both problems, which will
	be used to derive conditions concerning the switching points of the control function. A non-local optimality condition 
	treating back-and-forth switches will be formulated.

	For the numerical solution, we propose a proximal-gradient method. The emerging discretized subproblems will be solved 
	by employing Bellman's optimality principle, leading to an algorithm which is polynomial in the mesh size and in the admissible control 
	levels.
	An adaption of this algorithm can be used to handle subproblems of the trust-region method proposed in \cite{LeyfferManns2021}.
	Finally, we demonstrate computational results.
 \end{abstract}

 \begin{keywords}
	integer optimal control problem,
	total variation regularization,
	switching point optimization,
	proximal-gradient method,
	trust-region method
 \end{keywords}

\begin{msc}
\mscLink{49K30},
\mscLink{49L20}
\mscLink{49M37},
\mscLink{90C10}
\end{msc}

\section{Introduction}
\label{sec:intro}
We investigate the infinite-dimensional mixed-integer optimization problem
\begin{equation}
	\label{eq:prob}
	\tag{P}
	\begin{aligned}
		\text{Minimize} \quad & F(u) + \beta \TV(u) \\
		\text{such that} \quad & u(t) \in \set{ \nu_1, \ldots, \nu_d } \text{ for a.a.\ } t \in (0,T).
	\end{aligned}
\end{equation}
Here, the admissible control values satisfy $\set{ \nu_1, \ldots, \nu_d } \subset \Z$
with $\nu_1 < \nu_2 < \ldots < \nu_d$,
and
$\TV(u)$ is the total variation of the function $u$, see \cref{sec:prelim}.
The first part of the objective is kept rather general
and might contain, e.g., the solution operator of a differential equation.
Therefore,
\eqref{eq:prob}
covers a large class of mixed-integer optimal control problems
and these have an abundance of applications.
We refer to \cite{LeyfferManns2021}, \cite{SeverittManns2022} and the references therein. 

In \cite{LeyfferManns2021}, problems of the form \eqref{eq:prob} have been investigated and a trust-region algorithm has been proposed, with subproblems being modeled as linear integer problems. 
Here, we will extend some of the gained results.
For further investigation on mixed integer optimal control problems, see e.g. \cite{HanteSager2013}, 
\cite{BestehornHansknechtKirchesManns2020}, \cite{KirchesMannsUlbrich2021} and \cite{SagerZeile2021} using an 
approach based on the combinatorial integral approximation decomposition.

At this point, we would like to mention
that the total variation term in \eqref{eq:prob}
ensures the existence of minimizers under rather mild assumptions on $F$.
To be precise,
it suffices to assume that
$F \colon L^1(0,T) \to \R$ is lower semicontinuous and bounded from below,
see
\cite[Proposition~2.3]{LeyfferManns2021}
and the short argument after \cref{thm:props} below.
Since the total variation term penalizes the number (and height)
of the switches of the control function $u$,
it is also desirable from an application point of view.

The aim of this paper is threefold.
After recalling some properties of the total variation in \cref{sec:prelim},
we address optimality conditions for \eqref{eq:prob}
in \cref{sec:optimality_conditions}.
In particular,
we verify local optimality condition of first and second order
(\cref{thm:no_gap_SOC})
and we also formulate some non-local optimality conditions
(\cref{subsec:non-local_optimality})
in the spirit of the classical mode-insertion as in \cite[Section IV]{EgerstedtWardiAxelsson2006}.
Second,
we propose a proximal-gradient method for the solution of \eqref{eq:prob}
in \cref{sect:Proximal_Gradient}.
Third,
we give a fast algorithm for the solution
of the proximal-gradient subproblem (\cref{sect:prox_subproblems})
as well as for the subproblem arising in the trust-region method
proposed in \cite{LeyfferManns2021} (\cref{sec:TR}).
Finally,
we illustrate our findings by some numerical experiments in
\cref{sec:numerics}.

\section{The total variation functional}
\label{sec:prelim}
In this section,
we recall the definition of the total variation functional
$\TV \colon L^1(0,T) \to [0,\infty]$
and give some basic properties.

\begin{definition}
	\label{not:TV}
	Let $u\in L^1(0,T)$ and $a,b\in[0,T]$, $a<b$. Then,
	\[\TV(u;(a,b)):=\sup\set*{\int_a^b u \varphi' \, \d t\given \varphi\in C^1_c(a,b),~\norm{\varphi}_{L^\infty(a,b)} \le 1}.\]
	Furthermore, we write $\TV(u) := \TV(u;(0,T))$.
\end{definition}
The space of functions with bounded variation $\BV(0,T)$ is therefore defined as 
the set of all $u\in L^1(0,T)$ with $\TV(u)<\infty$, equipped with the norm
\[\norm{u}_{\BV(0,T)} = \norm{u}_{L^1(0,T)} + \TV(u).\]
Since both $F$ and $\TV$ are defined on $L^1(0,T)$,
we will ignore null sets in the following.

For the next sections, some properties of $\BV(0,T)$ are needed. 
\begin{theorem} 
	\label{thm:props}
	The space $\BV(0,T)$ and the functional $\TV$ have the following properties.
	\begin{enumerate}
		\item\label{thm:props:0}
			The space $\BV(0,T)$ is (isometric isomorphic to) the dual space of a separable Banach space.
		\item\label{thm:props:1} For a sequence $\seq{u_k}_{k\in\N}\subset \BV(0,T)$, we have $u_k\weaklystar 
		u$ in $\BV(0,T)$ if and only if $u_k\to u$ in $L^1(0,T)$ and 
		$\seq{u_k}_{k\in\N}$ is bounded in $\BV(0,T)$.
		\item\label{thm:props:2} $\BV(0,T)$ is continuously embedded in $L^\infty(0,T)$ and compactly embedded in $L^p(0,T)$ for all $p\in[1,\infty)$.
		\item\label{thm:props:3} When $u_k \weaklystar u$ in $\BV(0,T)$, we have $u_k\to u$ in $L^p(0,T)$ for all $p\in[1,\infty)$.
		\item\label{thm:props:4} If $\seq{u_k}_{k\in\N}$ is bounded in $\BV(0,T)$, there exists a weak-$\star$ accumulation point of $\seq{u_k}$.
		\item\label{thm:props:5}
			The functional $\TV$ is lower semicontinuous on $L^1(0,T)$,
			i.e., $u_k \to u$ in $L^1(0,T)$ implies $\TV(u) \le \liminf_{k \to \infty} \TV(u_k)$.
	\end{enumerate}
\end{theorem}
\begin{proof}
	For \ref{thm:props:0}, \ref{thm:props:1} and \ref{thm:props:2}, 
	see \cite[Remark 3.12, Proposition 3.13 and Corollary 3.49]{Ambrosio2000}.
	To prove \ref{thm:props:3}, we note that $\norm{u_k-u}_{L^1(0,T)}\to 0$ as 
	well as the boundedness of $\seq{u_k-u}_{k\in\N}$ in $\BV(0,T)$ follows from 
	\ref{thm:props:1}. Considering \ref{thm:props:2}, an interpolation 
	inequality yields that
	\[\norm{u_k-u}_{L^p(0,T)} \leq \norm{u_k-u}_{L^1(0,T)}^{1/p}\norm{u_k-u}_{L^\infty(0,T)}^{1-1/p}\to 0.\]

	Assertion \ref{thm:props:4} is a direct consequence of \ref{thm:props:0}.

	In order to prove \ref{thm:props:5},
	we take a subsequence with $\liminf_{k \to \infty} \TV(u_k) = \lim_{l \to \infty} \TV(u_{k_l})$.
	For an arbitrary $\varphi \in C_c^1(0,T)$ with $\norm{\varphi}_{L^\infty(0,T)} \le 1$,
	we have
	\begin{equation*}
		\int_0^T u \varphi' \, \dt
		=
		\lim_{l \to \infty}
		\int_0^T u_{k_l} \varphi' \, \dt
		\le
		\lim_{l \to \infty}
		\TV(u_{k_l})
		=
		\liminf_{k \to \infty} \TV(u_k)
		.
	\end{equation*}
	Taking the supremum over all these $\varphi$,
	we get the desired inequality.
\end{proof}
We define the set of admissible controls via
\[\Uad := \set*{u\in L^1(0,T)\given u(t)\in\set{\nu_1,\dots,\nu_d} \text{ for a.e.\ } t \in (0,T) }.\]
The
existence of a solution can be shown by standard arguments:
A minimizing sequence
$\seq{u_k}_{k\in\N}\subset\Uad$ is bounded in 
$L^1(0,T)$ by $T\max(\abs{\nu_1},\abs{\nu_d})$, while the boundedness of 
$\TV(u_k)$ follows from the existence of a lower bound for $F$.
Using \cref{thm:props} \ref{thm:props:4}, the existence of a weak-$\star$ 
convergent subsequence $\seq{u_{k_l}}_{l\in\N}$ with $u_{k_l}\weaklystar \bar 
u\in\BV(0,T)$ can be derived. Considering \cref{thm:props} \ref{thm:props:1}, 
we see that $u_{k_l}\to \bar u$ in $L^1(0,T)$. Thus, there is another 
subsequence $\seq{u_m}_{m\in\N}\subset \seq{u_{k_l}}_{l\in\N}$ with $u_m(t)\to 
\bar u(t)$ for a.e.\ $t\in(0,T)$. It follows that 
$\bar u(t)\in\set{\nu_1,\dots,\nu_d}$ a.e.\ in $(0,T)$, hence $\bar u\in\Uad$.
Finally, the lower semicontinuity of $F$ and \cref{thm:props} \ref{thm:props:5}
yield the optimality of $\bar u$.

The following lemma will be needed in \cref{sec:optimality_conditions}.
\begin{lemma}
	\label{lem:TV_inequality}
	Let $u\in\BV(0,T)$ and real values $0\le t_1<\dots<t_n\le T$ be given. Then, we have
	\begin{equation}
		\label{eq:TV_est}
		\sum_{j=1}^{n-1} \TV(u;(t_{j},t_{j+1}))
		\leq
		\TV(u)
		.
	\end{equation}
\end{lemma}
\begin{proof}
	By definition,
	there exist sequences $\seq{\varphi_{j,k}}_{k\in\N}\subset C_c^1(t_j,t_{j+1})$ for all $j\in\set{1,\dots,n-1}$ such that 
	$\norm{\varphi_{j,k}}_\infty\leq 1$ for all $k\in\N$ and
	\[
		\TV(u;(t_{j},t_{j+1}))
		=
		\lim_{k\to\infty}\int_{t_j}^{t_{j+1}} u
		\varphi_{j,k}'
		\, \d t
		.
	\]
	Then, we have
	\begin{equation*}
		\sum_{j=1}^{n-1} \TV(u;(t_{j},t_{j+1}))
		=
		\lim_{k\to\infty}\sum_{j=1}^{n-1}\int_{t_j}^{t_{j+1}} u \varphi_{j,k}'\, \d t
		=
		\lim_{k \to \infty} \int_0^T u \varphi_k' \, \dt
		\leq
		\TV(u),
	\end{equation*}
	where $\varphi_k \in C_c^1(0,T)$ is given by
	\begin{equation*}
		\varphi_k(t)
		=
		\begin{cases}
			\varphi_{j,k}(t) & \text{if } t \in (t_j, t_{j+1}) \text{ for some } j \in \set{1,\ldots,n-1} \\
			0 & \text{else}.
		\end{cases}
	\end{equation*}
\end{proof}
Note that we will not have equality in \eqref{eq:TV_est},
even in case $t_1 = 0$, $t_n = T$,
since jumps at the points $t_2, \ldots, t_{n-1}$
are ignored by the left-hand side of \eqref{eq:TV_est}.

\section{Optimality conditions}
\label{sec:optimality_conditions}
In this section, we are discussing optimality conditions
for \eqref{eq:prob}.
First, we address a switching-point reformulation in \cref{subsec:switching_points}.
This can be used to derive local optimality conditions of first and second order in \cref{subsec:local_opt_con}.
Afterwards,
we consider non-local optimality conditions in \cref{subsec:non-local_optimality}.

\subsection{Switching point reformulation}
\label{subsec:switching_points}
Let $n \in \N$, $t \in \R^{n-1}$, $a \in \R^n$
be given such that
$t_i \le t_{i+1}$ for all $i = 0, \ldots, n-1$,
with the convention $t_0 = 0$, $t_n = T$.
We define the function $v^{t,a} \in L^1(0,T) \cap \BV(0,T)$ via
\begin{equation*}
	v^{t,a}
	:=
	\sum_{j = 1}^n a_j \chi_{[t_{j-1}, t_{j})},
\end{equation*}
where we again use $t_0 = 0$ and $t_n = T$.

In \cite[Corollary~4.4]{LeyfferManns2021}
it is shown that each $u \in \Uad \cap \BV(0,T)$
has a (unique) representation $u = v^{t,a}$, $t \in \R^{n-1}$, $a \in \R^n$,
where $n$ is chosen as small as possible.
We give a different representation.
\begin{lemma}
	\label{lem:structure}
	Let $u \in \BV(0,T)$ be feasible for \eqref{eq:prob}.
	Then, there exists a unique $\hat{n}\in \N$
	and unique
	$\hat{t}_1, \ldots, \hat{t}_{\hat{n}-1} \in [0,T]$,
	$\hat{\kappa}_1, \ldots, \hat{\kappa}_{\hat{n}} \in \set{1, \ldots, d}$ 
	satisfying
	\begin{enumerate}
		\item\label{lem:structure:1}
		$0  < \hat{t}_1 \le \hat{t}_2 \le \ldots \le \hat{t}_{\hat{n}-1} < T$,
		\item\label{lem:structure:2}
		$\abs{\hat{\kappa}_j - \hat{\kappa}_{j+1}} = 1$
		for all $j = 1,\ldots, \hat{n}-1$,
		\item\label{lem:structure:3}
		$u = v^{\hat t, \hat a}$,
		where $\hat{a}_j = \nu_{\hat{\kappa}_j}$ for all $j=1,\dots,\hat{n}$,
		\item\label{lem:structure:4}
		if $\hat{t}_j = \hat{t}_{j+2}$, then $\hat{a}_j\neq \hat{a}_{j+2}$ for $j\in\{1,\dots,\hat{n}-3\}$.
	\end{enumerate}
\end{lemma}

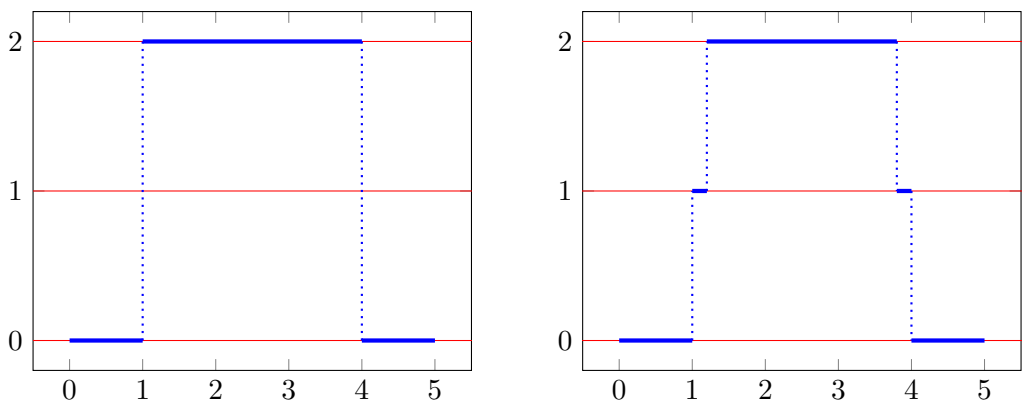
\begin{figure}[htp]
	\centering
	\begin{tikzpicture}
		\begin{axis}[width=.5\textwidth,xmin=-0.5,xmax=5.5, xtick ={0,1,2,3,4,5}]
			\addplot[red, mark={}] coordinates{ (-1,0) (6,0) };
			\addplot[red, mark={}] coordinates{ (-1,1) (6,1) };
			\addplot[red, mark={}] coordinates{ (-1,2) (6,2) };
			\addplot[blue, ultra thick, mark={}] coordinates{
					(0,0) (1,0)
					
					(1,2) (4,2)
					
					(4,0) (5,0)
				};

			\addplot[blue, thick, dotted, mark={}] coordinates{
				(1,0) (1,2)
				
				(4,2) (4,0)
			};
		\end{axis}
	\end{tikzpicture}%
	\hspace{1cm}%
	\begin{tikzpicture}%
		\begin{axis}[width=.5\textwidth,xmin=-0.5,xmax=5.5,xtick ={0,1,2,3,4,5}]
			\addplot[red, mark={}] coordinates{ (-1,0) (6,0) };
			\addplot[red, mark={}] coordinates{ (-1,1) (6,1) };
			\addplot[red, mark={}] coordinates{ (-1,2) (6,2) };
			\addplot[blue, ultra thick, mark={}] coordinates{
					(0,0) (1,0)
					
					(1,1) (1.2,1)
					
					(1.2,2) (3.8,2)
					
					(3.8,1) (4,1)
					
					(4,0) (5,0)
			};
			\addplot[blue, thick, dotted, mark={}] coordinates{
				(1,0) (1,1)	
				
				(1.2,1) (1.2,2)	
				
				(3.8,2) (3.8,1)	
						
				(4,1) (4,0)
			};
		\end{axis}
	\end{tikzpicture}%
	\caption{%
		The functions $u$ (left) and $\tilde u$ (right), see \cref{ex:Notation}.
		The control levels $\nu_1, \nu_2, \nu_3$ are visualized in red.%
	}
	\label{fig:hill_func}
\end{figure}

Before giving the proof,
we will explain the meaning of the conditions
\ref{lem:structure:1}--\ref{lem:structure:4}.
Using conditions \ref{lem:structure:1} and \ref{lem:structure:3}, we can identify $u$ with a piecewise constant function
with the switching points $\hat{t}_j$, $j\in\set{1,\dots,\hat{n}-1}$. In contrast to the
representation in \cite[Proposition~4.4]{LeyfferManns2021}, we also allow
equality of time steps. With \ref{lem:structure:2}, the equality of two or more $\hat{t}_j$ is needed
when $u$ is increasing or decreasing by more than one level.
Finally, \ref{lem:structure:4} prevents unnecessary and repetitive switching between two levels
at the same time instance.
To illustrate the difference to \cite[Proposition~4.4]{LeyfferManns2021}, we consider the following example.
\begin{example}
	\label{ex:Notation}
	We consider the situation with $d = 3$ control levels
	and $\set{\nu_1, \nu_2, \nu_3} = \set{0,1,2}$.
	For $T=5$, the function $u$ illustrated in \cref{fig:hill_func} (left) can be represented as
	$u = v^{t,a}$ or $u = v^{\hat t, \hat a}$ with
	\begin{align*}
		(t_1, t_2) &= (1, 4),
		&
		(\hat t_1, \hat t_2, \hat t_3, \hat t_4) &= (1, 1, 4, 4),
		\\
		(a_1, a_2, a_3) &= (0, 2, 0),
		&
		(\hat a_1, \hat a_2, \hat a_3, \hat a_4, \hat a_5) &= (0, 1, 2, 1, 0).
	\end{align*}
	While the second representation of $u$ seems to be overcomplicated,
	the function $\tilde u$ from \cref{fig:hill_func} (right) can be represented using the
	$\hat{a}_j$ defined above by simply adapting the time steps.
	Indeed, we have
	$\tilde u = v^{\tilde t, \hat a}$
	with $(\tilde t_1, \tilde t_2, \tilde t_3, \tilde t_4) = (1,1.2,3.8,4)$.
	Note that $\tilde u$ can be interpreted as a perturbation of the original function $u$.
	This is not possible by using the first representation of $u$,
	since this representation does not include the control level $\nu_2 = 1$.
\end{example}

\begin{proof}[Proof of \cref{lem:structure}]
	In \cite[Proposition 4.4]{LeyfferManns2021},
	the existence of $n\in\N$ and $t_1,\dots,t_{n-1}$, $\kappa_1,\dots,\kappa_n$ satisfying
	$0<t_1<\dots<t_{n-1} < T$,
	$\kappa_j\neq \kappa_{j+1}$ for $j\in\{1,\dots,n-1\}$
	and \ref{lem:structure:3} has been proven.
	To fulfil \ref{lem:structure:1}--\ref{lem:structure:4},
	we can construct the time steps $\hat{t}_1,\dots,\hat{t}_{\hat{n}-1}$ by appending for every $j\in\set{1,\dots,n-1}$ with $\abs{\kappa_{j+1}-\kappa_j}>1$
	in total $\abs{\kappa_{j+1}-\kappa_j}-1$ new time steps equal to $t_j$
	such that \ref{lem:structure:2} is accomplished.
	Notice that this  method implies \ref{lem:structure:4}
	since we added the minimum number of time steps needed to ascend or descend from $\kappa_j$ to $\kappa_{j+1}$,
	while \ref{lem:structure:3}  is still valid, considering the characteristic function of the empty set equals zero.

	The uniqueness of $\hat n, \hat t, \hat a$ is easy to check.
\end{proof}
In what follows, we associate with a given function
the representations from \cite[Proposition~4.4]{LeyfferManns2021}
and from \cref{lem:structure}.
\begin{notation}
	\label{not:u_and_switching_time}
	Let $u \in \BV(0,T)$ be feasible for \eqref{eq:prob}.
	First, we use \cite[Proposition~4.4]{LeyfferManns2021},
	to get the representation
	$u = v^{t,a}$
	and $a_j = \nu_{\kappa_j}$ for some $\kappa_j \in \set{1,\ldots,d}$.
	Here, the value $n \in \N$ is as small as possible,
	thus, we refer to
	$u = v^{t,a}$
	as the \emph{minimal representation} of $u$.

	Second, we use \cref{lem:structure} to get the representation
	$u = v^{\hat t, \hat a}$
	and $\hat a_j = \nu_{\hat \kappa_j}$ for some $\hat \kappa_j \in 
	\set{1,\ldots,d}$.
	Here, the changes
	$\abs{\hat a_j - \hat a_{j-1}}$
	(or, equivalently, $\abs{\hat \kappa_j - \hat \kappa_{j-1}}$)
	are as small as possible,
	cf.\ \itemref{lem:structure:2}.
	This means that the jumps are fully resolved
	and, therefore,
	we refer to
	$u = v^{\hat t, \hat a}$
	as the \emph{full representation} of $u$.

	Finally, we define the index sets
	(associated with the minimal representation)
	\begin{equation*}
		J^+ = \set{j \in \set{1,\ldots,n} \given \kappa_{j+1} > \kappa_j + 1},\quad
		J^- = \set{j \in \set{1,\ldots,n} \given \kappa_{j+1} < \kappa_j - 1}.
	\end{equation*}
	The set $J^+$ ($J^-$) consists of exactly those indices $j$,
	for which there is an upwards (downwards) jump at $t = t_j$
	which skips over the control levels between $a_j$ and $a_{j+1}$.
\end{notation}
Using the full representation of a feasible function, the following can be proved.
\begin{lemma}
	\label{lem:similar_local_points}
	Let $u$ be a feasible point of \eqref{eq:prob} and $\hat t\in\R^{\hat n-1}$, $\hat a\in\R^{\hat n}$ be chosen such that $u = v^{\hat t, \hat a}$ is the full representation of $u$.
	Then, there exist $\hat\varepsilon,c,C > 0$ such that for all feasible points $w$ of \eqref{eq:prob}
	with $\TV(w) \le \TV(u)$ and $\norm{w - u}_{L^1(0,T)} \le \hat\varepsilon$
	there exists $\hat s \in \R^{\hat{n} - 1}$
	with
	$0 < \hat s_1 \le \ldots \le \hat s_{\hat{n}-1} < T$
	and
	$w = v^{\hat s, \hat a}$.
	Furthermore, the estimate
	\begin{equation*}
		c \norm{\hat s - \hat t}_{\R^{\hat{n}-1}}
		\le
		\norm{w - u}_{L^1(0,T)}
		\le
		C \norm{\hat s - \hat t}_{\R^{\hat{n}-1}}
	\end{equation*}
	holds.
\end{lemma}
\begin{proof}
	We set
	\[
	\hat\varepsilon := \tfrac 12 \min\set*{\hat{t}_{j+1}-\hat{t}_j\given \hat{t}_j\neq \hat{t}_{j+1},~ j\in\{0,\dots,\hat n-1\}}.
	\]
	Let $w$ be a feasible point of \eqref{eq:prob} with $\TV(w) \le \TV(u)$
	and $\norm{w - u}_{L^1(0,T)} \le \hat\varepsilon$.
	Thus,
	\begin{equation*}
		\hat\varepsilon
		\ge
		\norm{w - u}_{L^1(0,T)}
		\ge
		\lambda(\set{t \in (0,T) \given u(t) \ne w(t)}),
	\end{equation*}
	where $\lambda$ is the Lebesgue measure.
	Since $u$ and $w$ are piecewise constant,
	there exists a nonempty interval 
	$(\hat\alpha_j,\hat\beta_j)\subset [\hat t_{j}, \hat t_{j+1}]$ for every 
	$j\in\{0,\dots,\hat n-1\}$ with $\hat t_{j}\neq \hat t_{j+1}$ on which $w = u = \hat a_{j}$.
	The same is 
	true when considering the minimal representation $v^{t,a}$ with $t\in\R^{n-1}$, $a\in\R^n$ of $u$, where we 
	get the existence of such an interval in $[t_{j},t_{j+1}]$ for every $j\in\{0,\dots,n-1\}$ on which $w = u = a_j$.

	Let $w=v^{s,\tilde a}$ be the minimal representation of $w$ with $s\in\R 
	^{m-1}$, $\tilde a = (\tilde a_1,\dots,\tilde a_m)$. Since there is an open 
	subinterval $(\alpha_j,\beta_j)$ of $[t_{j},t_{j+1}]$ with $w = a_j$, we can 
	define the midpoint $\tilde t_j$ of this interval for every $j\in\{0,\dots 
	,n-1\}$. By defining $\varphi\in C_c^1(\tilde t_j,\tilde t_{j+1})$ as a 
	continuous function with $\varphi(t) = -\sgn(a_{j+1}-a_j)$ for 
	$t\in(\beta_j,\alpha_{j+1})$,
	we can see that
	\begin{equation*}
		\begin{aligned}
			\TV(w; (\tilde t_{j},\tilde t_{j+1}))&\geq \int_{\tilde t_j}^{\tilde t_{j+1}}u\varphi'\, \d t = a_j\int_{\tilde t_j}^{\beta_j}\varphi'\, \d t + a_{j+1}\int_{\alpha_{j+1}}^{\tilde t_{j+1}}\varphi'\, \d t\\
			&= -\sgn(a_{j+1}-a_j)(a_j - a_{j+1}) = \abs{a_{j+1}-a_j}.
		\end{aligned}
	\end{equation*}
	Then, using \cref{lem:TV_inequality}, it follows that that
	\[
		\TV(u) \ge \TV(w) \ge \sum_{j=0}^{n-1} \TV(w; (\tilde t_{j},\tilde t_{j+1}))\geq 
	\sum_{j=0}^{n-1} \abs{a_{j+1} - a_j} = \TV(u).
	\]
	Thus, equality holds. In particular, we have 
	\[
	\TV(w; (\tilde t_{j}, \tilde t_{j+1})) = \abs{a_{j+1}-a_j},
	\]
	implying that $w$ can only ascend or descend from $a_{j}$ to $a_{j+1}$ in 
	$(\tilde t_{j},\tilde t_{j+1})$. Translating this behaviour in the full representation, we see that for every $j\in\{0,\dots,\hat n-1\}$ with $\hat t_{j}\neq \hat t_{j+1}$, $w$ has to switch to every value between $\hat a_j$ and $\hat a_{j+1}$ exactly once in $(\hat\alpha_j,\hat\beta_{j+1})$. We conclude that the full representation of $w$ is given by $v^{\hat s,\hat a}$ for an $\hat s\in\R^{\hat n-1}$.
	
	Now, observe that
	\begin{equation}
		\label{eq:wmu}
		w - u
		=
		v^{\hat s,\hat a} - v^{\hat t, \hat a}
		=
		\sum_{j = 1}^{\hat n-1} \mu_j \sgn(\tau_j) \chi_{I_j(\tau)},
	\end{equation}
	with $\tau_j = \hat s_j - \hat t_j$,
	$\mu_j = \hat a_{j+1} - \hat a_j$
	and
	\begin{equation*}
		I_j(\tau) =
		\begin{cases}
			(\hat t_j, \hat s_j) & \text{if } \tau_j > 0, \\
			\emptyset & \text{if } \tau_j = 0, \\
			(\hat s_j, \hat t_j) & \text{if } \tau_j < 0.
		\end{cases}
	\end{equation*}
	Note that, at every $t \in (0,T)$,
	all non-vanishing addends on the right-hand side of \eqref{eq:wmu}
	share the same sign.
	Thus,
	\begin{equation*}
		\norm{ w - u }_{L^1(0,T)}
		=
		\sum_{j = 1}^{\hat n-1} \norm{\mu_j \sgn(\tau_j) \chi_{I_j(\tau)}}_{L^1(0,T)}
		=
		\sum_{j = 1}^{\hat n} \abs{\hat s_j-\hat t_j} \abs{\hat a_j - \hat a_{j-1}}
		.
	\end{equation*}
	Since $\hat a_j \in \{\nu_1,\dots,\nu_d\}\subset\Z$ $\forall j\in\set{1,\dots,\hat n}$, we conclude
	\[
		\sum_{j = 1}^{\hat n} \abs{\hat s_j-\hat t_j}
		\leq
		\norm{w-u}_{L^1(0,T)}
		\leq
		(\nu_d-\nu_1)
		\sum_{j = 1}^{\hat n} \abs{\hat s_j-\hat t_j}
	\]
 	from which, using the equivalence of all norms in $\R^{\hat n-1}$, the statement follows.
\end{proof}
Now, we want to derive local optimality conditions for \eqref{eq:prob}
via
reformulation as a switching point optimization problem similar to \cite[Section~4.2]{LeyfferManns2021}.
Given $n \in \N$ and $a \in \R^n$, we consider the problem
\begin{equation}
	\label{eq:ST}
	\tag{ST($n,a$)}
	\begin{aligned}
		\text{Minimize} \quad & F(v^{t,a}) \\
		\text{with respect to} \quad & t \in \R^{n-1},\\
		\text{such that} \quad
		& 0 \le t_1 \le \dots \le t_{n-1} \le T.
	\end{aligned}
\end{equation}
Note that \eqref{eq:ST}
depends on the chosen values of $n \in \N$ and $a \in \R^{n}$.
We mention that we also utilize
\STfull,
where we use the data $(\hat n, \hat a)$ from the full representation of $u$.
The main advantage of using the full representation
is the upcoming theorem
showing that local optimality of $u$ for \eqref{eq:prob}
is equivalent to local optimality of $\hat t$ for \STfull.

\begin{theorem}
	\label{thm:equivalence}
	Let $u \in \BV(0,T)$ be feasible for \eqref{eq:prob}
	and consider the data $(\hat n, \hat a, \hat t)$ of its full representation.
	Then, $u$ is locally optimal for \eqref{eq:prob} in $L^1(0,T)$
	if and only if $\hat t$ is locally optimal for \STfull.
	Moreover, $u$ satisfies a local quadratic growth condition for \eqref{eq:prob}
	in $L^1(0,T)$
	if and only if a local quadratic growth condition is valid for \STfull\ at $\hat t$.
	To be precise, the existence of constants $\varepsilon,\eta > 0$
	with
	\begin{equation}
		\label{eq:quad_growth_P}
		F(w) + \beta \TV(w)
		\ge
		F(u) + \beta \TV(u) + \frac\eta2 \norm{w - u}_{L^1(0,T)}^2
		\quad
		\forall w \in \Uad, \norm{w - u}_{L^1(0,T)} \le \varepsilon
	\end{equation}
	is equivalent to the existence of constants $\tilde\varepsilon, \tilde\eta > 0$
	with
	\begin{equation}
		\label{eq:quad_growth_STfull}
		F(v^{\hat s,\hat a})
		\ge
		F(u) + \frac{\tilde\eta}2 \norm{\hat s - \hat t}_{\R^{\hat n-1}}^2
		\qquad
		\forall \hat s \in \FF,
		\norm{\hat s - \hat t}_{\R^{\hat n - 1}} \le \tilde\varepsilon
		,
	\end{equation}
	where
	\begin{equation*}
		\FF := \set{ \hat s \in \R^{\hat n - 1} \given 0 \le \hat s_1 \le \dots \le \hat s_{\hat n - 1} \le T }
	\end{equation*}
	is the feasible set of \STfull.
\end{theorem}
\begin{proof}
	We suppose that $u = v^{\hat t, \hat a}$ satisfies \eqref{eq:quad_growth_P} with $\varepsilon > 0$
	and $\eta \ge 0$.
	Note that $\eta = 0$ corresponds to local optimality of $u$,
	whereas $\eta > 0$ describes a quadratic growth condition.
	Similar to the proof of
	\cite[Lemma~4.12]{LeyfferManns2021}, we define 
	\[h:=\min\set*{\tfrac 12 \min\set{\hat t_{i+1}-\hat t_i\given \hat t_{i+1}\neq \hat t_i,~i\in\{0,\dots,n-1\}},\frac{\varepsilon}{\hat n \parens{\nu_d-\nu_1}}} > 0\]
	and choose $\tilde\varepsilon \in (0,h)$.
	Then, for every $\hat s\in \FF \cap B_{\tilde\varepsilon}(\hat t)$
	we have by construction $\TV(u) = \TV(v^{\hat s,\hat a})$ as well as 
	$\norm{u - v^{\hat s,\hat a}}_{L^1(0,T)}
	= \norm{v^{\hat t,\hat a} - v^{\hat s,\hat a}}_{L^1(0,T)}\leq \abs{\nu_d-\nu_1}\hat n 
	h\leq \varepsilon$, thus, 
	\begin{equation*}
		F(v^{\hat s,\hat a})
		\geq
		F(u)
		+
		\frac\eta2 \norm{ u - v^{\hat s, \hat a}}_{L^1(0,T)}^2
		\geq
		F(u)
		+
		\frac{\eta c^2}2 \norm{ \hat s - \hat t }_{\R^{\hat n-1}}^2
		,
	\end{equation*}
	where we used \cref{lem:similar_local_points}.
	Thus, we arrive at \eqref{eq:quad_growth_STfull} with $\tilde\eta = \eta c^2$.
	Note that $\tilde\eta > 0$ if $\eta > 0$.
	This shows local optimality of $\hat t$ if $\eta = 0$
	and the quadratic growth condition if $\eta > 0$.

	For the converse implications, we assume that $\hat t$
	satisfies \eqref{eq:quad_growth_STfull} with $\tilde\varepsilon > 0$ and $\tilde\eta \ge 0$.
	We define $\eta := \tilde\eta/C^2$ with $C > 0$ from \cref{lem:similar_local_points}.
	Note that the continuity of $w \mapsto \norm{w - u}^2_{L^1(0,T)}$ 
	implies that $w\mapsto F(w)-\tfrac\eta 2 \norm{w - u}^2_{L^1(0,T)}$ is lower 
	semicontinuous. Hence,
	\[\tilde M:=\set*{w\in L^1(0,T)\given F(w)-\frac\eta 2 \norm{w - u}^2_{L^1(0,T)}>F(u)-\beta}\]
	is open, and due to $u\in\tilde M$ there exists $\bar\varepsilon>0$ with $B_{\bar\varepsilon}(u)\subset \tilde M$.
	We define
	$\varepsilon := \min\set{c\tilde\varepsilon,\hat\varepsilon,\bar\varepsilon}$ 
	with $\hat\varepsilon$ and $c$ given by \cref{lem:similar_local_points}.
	Let $w\in\Uad$ with $\norm{v^{\hat t,\hat a}-w}_{L^1(0,T)}\leq \varepsilon$,
	be given.
	In case
	$\TV(w)\leq \TV(v^{\hat t,\hat a})$ there exists $\hat s \in \FF$ with 
	$w=v^{\hat s,\hat a}$
	and we have
	$\norm{\hat s - \hat t}_{\R^{\hat n - 1}} \le \frac1c \norm{w - u}_{L^1(0,T)} \le \tilde \varepsilon$.
	Thus, \eqref{eq:quad_growth_STfull} and \cref{lem:similar_local_points} yield
	\begin{equation*}
		\begin{aligned}
			F(w)+\beta\TV(w)&\geq F(u)+\beta\TV(u) + \frac{\tilde\eta}{2}\norm{\hat s - \hat t}_{\R^{\hat n-1}}^2\\
			&\geq F(u)+\beta\TV(u) + \frac{\eta}{2}\norm{w - u}_{L^1(0,T)}^2
			.
		\end{aligned}
	\end{equation*}
	In the other case $\TV(w)>\TV(u)$, we use $w \in \tilde M$ to obtain
	\begin{equation*}
		\begin{aligned}
			F(w)+\beta\TV(w)&>F(u)-\beta + \frac\eta 2 \norm{w - u}^2_{L^1(0,T)} + \beta(\TV(u)+1) \\
			&= F(u)+\beta\TV(u) + \frac\eta 2 \norm{w - u}^2_{L^1(0,T)}
			.
		\end{aligned}
	\end{equation*}
	Hence, we have shown \eqref{eq:quad_growth_P}
	with $\eta \ge 0$.
\end{proof}
Note that equivalence of the local optimalities will not hold in general
if we are using the minimal representation.

\subsection{Local optimality conditions for \texorpdfstring{\eqref{eq:prob}}{(P)}}
\label{subsec:local_opt_con}
In this section, we derive optimality conditions for \eqref{eq:prob}
via the (equivalent) problem \STfull.
To this end,
we are going to discuss optimality conditions for the problem
\eqref{eq:ST} and
these findings will also be applied to \STfull.
Since \eqref{eq:ST}
is a standard finite-dimensional optimization problem,
optimality conditions involving
first and second order derivatives of the objective of \eqref{eq:ST}
(w.r.t.\ $t$)
can be formulated.
Thus, we are going to investigate these derivatives.

In the upcoming theorem,
we need some regularity of $F$.
First, we assume
that
$F \colon L^1(0,T) \to \R$
is twice Fréchet differentiable.
This yields the second-order Taylor expansion
\begin{equation*}
	F(v)
	=
	F(u) + F'(u)( v - u) + \frac12 F''(u) ( v - u)^2
	+
	\oo\parens*{\norm{v - u}_{L^1(0,T)}^2}
	\quad\text{as } \norm{v - u}_{L^1(0,T)} \to 0,
\end{equation*}
see
\cite[Theorem~5.6.3]{Cartan1967}.
Here, $F'(u) \colon L^1(0,T) \to \R$ and $F''(u) \colon L^1(0,T) \times L^1(0,T) \to \R$
are the Fréchet derivatives of first and second order at $u$, respectively,
and $F''(u) (v-u)^2$ is short for $F''(u)[v-u,v-u]$.
We investigate the structure of the derivatives.
The first order derivative $F'(u)$
belongs to the dual space of $L^1(0,T)$,
which will be identified with $L^\infty(0,T)$.
Thus, $F'(u)$ is identified with a function $\nabla F(u) \in L^\infty(0,T)$ and we will pose
regularity assumptions on this function.
Similarly, $F''(u)$ is a continuous bilinear form on $L^1(0,T)$.
It is well known that continuous bilinear forms on $L^1(0,T)$
can be identified with functions from $L^\infty( (0,T)^2)$.
In fact, this follows from the (isometric) identifications
\begin{align*}
	\mathfrak{Bil}(L^1(0,T), L^1(0,T))
	&\isometric
	(L^1(0,T) \mathbin{\otimes_\pi} L^1(0,T) )\dualspace
	=
	(L^1(0,T) \mathbin{\otimes_{\Delta_1}} L^1(0,T) )\dualspace
	\\&\isometric
	L^1( (0,T)^2 ) \dualspace
	\isometric
	L^\infty( (0,T)^2 )
	,
\end{align*}
see \cite[Sections~3 and 7]{DefantFloret1992} for the results and for the notation.
Thus, we will identify $F''(u)$ with a function $\nabla^2 F(u)$ from $L^\infty( (0,T)^2 )$
and the evaluation (given by the above identifications) is
\begin{equation*}
	F''(u) [v_1, v_2]
	=
	\int_0^1 \int_0^1 \nabla^2 F(u)(r,s) v_1(r) v_2(s) \, \d r \, \d s
	\qquad
	\forall v_1, v_2 \in L^1(0,T).
\end{equation*}
As for $\nabla F(u) \colon (0,T) \to \R$,
we are going to postulate
regularity assumptions on the function
$\nabla^2 F(u) \colon (0,T)^2 \to \R$.
Finally,
we mention that the symmetry of $F''(u)$,
see \cite[Theorem~5.1.1]{Cartan1967},
yields $\nabla^2 F(u)(t,s) = \nabla^2 F(u)(s,t)$
for a.a.\ $(s,t) \in (0,T)^2$.

\begin{theorem}
	\label{thm:second_derivatives}
	We consider fixed $n \in \N$, $a \in \R^{n}$.
	Let the vector $t \in \R^{n-1}$ be feasible for \eqref{eq:ST}
	and let $\tau \in \R^{n-1}$ be given such that $\tau_k \le \tau_{k+1}$
	whenever $t_k = t_{k+1}$ for all $k = 0,\ldots, n$
	with the convention $0 = t_0 = \tau_0 = \tau_n$ and $T = t_n$.
	Then, $t + \tau$ is feasible for \eqref{eq:ST}
	whenever $\norm{\tau}$ is small enough.
	Under the regularity assumptions
	that $F \colon L^1(0,T) \to \R$
	is twice Fréchet differentiable at $v^{t,a}$,
	$\nabla F(v^{t,a}) \in C^1([0,T])$
	and $\nabla^2 F(v^{t,a}) \in C([0,T]^2)$,
	we have the expansion
	\begin{align*}
		F(v^{t + \tau, a})
		&=
		F(v^{t,a})
		+
		\parens*{ \sum_{j = 1}^{n-1} \mu_j \nabla F(v^{t,a})(t_j) \tau_j }
		\\&\quad
		+
		\frac{1}{2} \parens*{
			\sum_{j = 1}^{n-1} \mu_j (\nabla F(v^{t,a}))'(t_j) \tau_j^2
			+
			\sum_{j,k = 1}^{n-1} \mu_j \mu_k \nabla^2 F(v^{t,a})(t_j, t_k) \tau_j \tau_k
		}
		\\&\quad
		+
		\oo( \norm{\tau}^2 )
		\qquad\text{as } \norm{\tau} \to 0
		.
	\end{align*}
	Here, $\mu_j = a_{j+1} - a_j$ is the jump height at $t_j$.
\end{theorem}
\begin{proof}
	The feasibility of $t + \tau$ for $\norm{\tau}$ small enough is clear.
	For brevity, we write $v^{t}$ and $v^{t+\tau}$
	instead of $v^{t,a}$ and $v^{t+\tau,a}$, respectively.
	By definition of $v^{t + \tau}$ and $v^{t}$, we have
	\begin{equation*}
		v^{t + \tau} - v^{t}
		=
		\sum_{j = 1}^{n-1} \mu_j \sgn(\tau_j) \chi_{I_j(\tau)},
	\end{equation*}
	with
	\begin{equation*}
		I_j(\tau) =
		\begin{cases}
			(t_j, t_j + \tau_j) & \text{if } \tau_j > 0, \\
			\emptyset & \text{if } \tau_j = 0, \\
			(t_j + \tau_j, t_j) & \text{if } \tau_j < 0.
		\end{cases}
	\end{equation*}
	Note that $\norm{ v^{t + \tau} - v^{t}}_{L^1(0,T)} \to 0$ as $\norm{\tau} \to 0$.
	Since $F$ is assumed to be twice Fréchet differentiable on $L^1(0,T)$, we get the expansion
	\begin{equation*}
		F(v^{t + \tau})
		=
		F(v^{t})
		+
		F'(v^{t}) (v^{t + \tau} - v^{t})
		+
		\frac12 F''(v^{t}) (v^{t + \tau} - v^{t})^2
		+
		\oo(\norm{v^{t + \tau} - v^{t}}_{L^1(0,T)}^2 )
		.
	\end{equation*}
	Note that
	$\oo(\norm{v^{t + \tau} - v^{t}}_{L^1(0,T)}^2 )$
	is already
	$\oo(\norm{\tau}^2)$
	due to
	$\norm{v^{t + \tau} - v^{t}}_{L^1(0,T)} \le C \norm{\tau}$.

	We study the terms on the right-hand side
	of the expansion
	by using the above representation of $v^{t + \tau} - v^{t}$.
	First, we have
	\begin{equation*}
		F'(v^{t}) (v^{t + \tau} - v^{t})
		=
		\sum_{j = 1}^{n-1} \mu_j \sgn(\tau_j) \int_{I_j(\tau)} \nabla F(v^{t})(s) \, \ds
		.
	\end{equation*}
	By using $\nabla F(v^{t})(s) = \nabla F(v^{t})(t_j) + (\nabla F(v^{t}))'(t_j) (s - t_j) + \oo( s - t_j )$,
	we find
	\begin{equation*}
		F'(v^{t}) (v^{t + \tau} - v^{t})
		=
		\sum_{j = 1}^{n-1} \mu_j
		\parens*{
			\nabla F(v^{t})(t_j) \tau_j
			+
			\frac12 (\nabla F(v^{t}))'(t_j) \tau_j^2
		}
		+ \oo( \tau_j^2 )
		.
	\end{equation*}
	Similarly,
	\begin{align*}
		F''(v^{t}) (v^{t + \tau} - v^{t})^2
		&=
		\sum_{j = 1}^{n-1} \sum_{k = 1}^{n-1}
		\mu_j \sgn(\tau_j) \mu_k \sgn(\tau_k)
		\int_{I_j(\tau)} \int_{I_k(\tau)} \nabla^2 F(v^{t})(r,s) \, \d r \, \d s
		\\
		&=
		\sum_{j = 1}^{n-1} \sum_{k = 1}^{n-1} \mu_j \mu_k \tau_j \tau_k \nabla^2 F(v^{t})(t_j, t_k)
		+
		\oo( \tau_j \tau_k )
		,
	\end{align*}
	where we used continuity of the function $\nabla^2 F(v^{t})$.
	This shows the claim.
\end{proof}
We note that the first order part of the expansion can be shown
by assuming first order Fréchet-differentiability of $F \colon L^1(0,T) \to \R$
at $v^{t,a}$
and continuity of $\nabla F(v^{t,a}) \colon [0,T] \to \R$.

\begin{lemma}
	\label{thm:optimality_conditions_ST}
	We consider fixed $n \in \N$, $a \in \R^{n}$.
	Let the vector $t \in \R^{n-1}$ be feasible for \eqref{eq:ST}
	and $t_0 = 0 < t_1$, $t_{n-1} < t_n = T$.
	We again use
	the jump heights
	$\mu_j := a_{j+1} - a_j$
	and define
	\begin{equation*}
		\TT :=
		\set[\big]{
			\tau \in \R^{n-1}
			\given
			\forall k \in \set{1,\ldots,n-2}
			:
			t_k = t_{k+1}
			\;\Rightarrow\;
			\tau_k \le \tau_{k+1}
		}
		.
	\end{equation*}
	We assume $\mu_j \ne 0$ for all $j = 1,\ldots,n-1$
	and
	we suppose that all jumps at $t_i$ go in the same direction,
	i.e.,
	$\sgn(\mu_i) = \sgn(\mu_j)$ for all $i,j \in \set{1,\ldots,n-1}$ with $t_i = t_j$.
	Further, we assume that $F$ satisfies the regularity assumptions
	of \cref{thm:second_derivatives}.
	If $t$ is a local minimizer of \eqref{eq:ST},
	then
	\begin{subequations}
		\label{eq:necessary}
		\begin{align}
			\label{eq:necessary_1}
			&\mathrlap{\forall j = 1,\ldots, n-1:}
			&
			\nabla F(v^{t,a})(t_j) &= 0
			,
			\\
			\label{eq:necessary_2}
			&\forall \tau \in \TT:
			&
			\sum_{j = 1}^{n-1} \mu_j (\nabla F(v^{t,a}))'(t_j) \tau_j^2
			+
			\sum_{j,k = 1}^{n-1} \mu_j \mu_k \nabla^2 F(v^{t,a})(t_j, t_k) \tau_j \tau_k
			&
			\ge0
			.
		\end{align}
	\end{subequations}
	On the other hand,
	if
	\begin{subequations}
		\label{eq:sufficient}
		\begin{align}
			\label{eq:sufficient_1}
			&\mathrlap{\forall j = 1,\ldots, n-1:}
			&
			\nabla F(v^{t,a})(t_j) &= 0
			,
			\\
			\label{eq:sufficient_2}
			&\forall \tau \in \TT\setminus\set{0}:
			&
			\sum_{j = 1}^{n-1} \mu_j (\nabla F(v^{t,a}))'(t_j) \tau_j^2
			+
			\sum_{j,k = 1}^{n-1} \mu_j \mu_k \nabla^2 F(v^{t,a})(t_j, t_k) \tau_j \tau_k
			&
			>0
		\end{align}
	\end{subequations}
	is satisfied,
	then $t$ is a local minimizer of \eqref{eq:ST}
	and a quadratic growth condition is satisfied.
\end{lemma}
The assumption $\mu_j \ne 0$ means that there is actually a jump at $t = t_j$
and the second assumption on $\mu$ corresponds to 
\itemref{lem:structure:4}.
\begin{proof}
	It is straightforward to verify that \eqref{eq:ST}
	satisfies the linear independence constraint qualification.
	This implies
	that $\TT$ coincides with the tangent cone of the feasible set at the point $t$,
	see \cite[Lemma~12.2]{NocedalWright2006}.
	Next, we are going to employ optimality conditions of first and second order.
	Note that there is a slight difficulty,
	since the objective of \eqref{eq:ST} is only defined on the
	feasible set, which is a closed set.
	However, we have proven a Taylor-like second order expansion
	in \cref{thm:second_derivatives}.
	By inspecting the proofs of
	\cite[Theorems~12.3, 12.5 and 12.6]{NocedalWright2006},
	we see that this is enough in order to get optimality conditions.

	To prove the necessary conditions,
	we assume that $t$ is locally optimal.
	The first-order optimality condition (\cite[Theorems~12.3]{NocedalWright2006}) reads
	\begin{equation*}
		\sum_{j = 1}^{n-1} \mu_j \nabla F(v^{t,a})(t_j) \tau_j
		\ge
		0
		\qquad\forall \tau \in \TT.
	\end{equation*}
	For any $j \in \set{1,\ldots,n}$,
	there exist $i,k \in \set{1,\ldots,n}$
	with $i \le j \le k$,
	\begin{equation*}
		t_{i-1} < t_i = t_j = t_k < t_{k+1}.
	\end{equation*}
	Then, the unit vectors
	$-e_i$ and $e_k$ belong to $\TT$ and this gives \eqref{eq:necessary_1}
	due to $\sgn(\mu_i) = \sgn(\mu_k)$.
	Since the derivative of the objective is zero,
	the critical cone used for second order conditions
	coincides with the tangent cone $\TT$ and the Lagrange multipliers are zero.
	The second-order necessary condition \cite[Theorem~12.5]{NocedalWright2006}
	delivers \eqref{eq:necessary_2}.

	The sufficiency of \eqref{eq:sufficient}
	follows with similar arguments from \cite[Theorem~12.6]{NocedalWright2006}.
\end{proof}

\begin{remark}
	\label{rem:different_SOC}
	If $u \in \BV(0,T)$ is feasible for \eqref{eq:prob}
	and has a switch across more than one level,
	i.e., if it switches from $\nu_i$ to $\nu_j$ with $\abs{i - j} > 1$,
	then
	the minimal representation $(t,a)$ and the full representation $(\hat t, \hat a)$
	deliver two different instances \eqref{eq:ST} and \STfull.
	It is easy to check that the first order part of \cref{thm:optimality_conditions_ST}
	gives the same conditions, namely
	$\nabla F(u)(t) = 0$ for all switching times $t \in (0,T)$.
	By means of an example, we check that the second order conditions differ.
	
	We consider the setting
	\begin{equation*}
		T = 2, \quad
		d = 3, \quad
		\nu_1 = 0, \quad
		\nu_2 = 1, \quad
		\nu_3 = 2
	\end{equation*}
	and the feasible point $u = \nu_1 \chi_{(0,1]} + \nu_3 \chi_{(1,2)}$,
	which has a jump from $\nu_1$ to $\nu_3$ at $t = 1$.
	The minimal representation of $u$ is given by
	\begin{equation*}
		n = 2, \quad
		a_1 = 0, \quad
		a_2 = 2, \quad
		t_1 = 1.
	\end{equation*}
	Consequently, $\TT = \R$
	and the condition \eqref{eq:necessary_2} reads
	\begin{equation}
		\label{eq:SOCLM}
		2 (\nabla F(u))'(1)
		+
		4 \nabla^2 F(u)(1,1)
		\ge 0.
	\end{equation}
	On the other hand,
	the full representation of $u$
	is given by
	\begin{equation*}
		\hat n = 3, \quad
		\hat a_1 = 0, \quad
		\hat a_2 = 1, \quad
		\hat a_3 = 2, \quad
		\hat t_1 = \hat t_2 = 1.
	\end{equation*}
	The application of \cref{thm:optimality_conditions_ST}
	to the full representation results in
	$\hat\TT = \set{\tau \in \R^2 \given \tau_1 \le \tau_2}$
	and
	\eqref{eq:necessary_2} is equivalent to
	\begin{equation}
		\label{eq:SOCWE}
		(\nabla F(u))'(1) ( \tau_1^2 + \tau_2^2)
		+
		\nabla^2 F(u)(1,1) (\tau_1 + \tau_2)^2
		\ge0
		\quad
		\forall \tau \in \hat\TT
		.
	\end{equation}
	For $\tau_1 = \tau_2$ this is exactly \eqref{eq:SOCLM},
	but since we can also choose $\tau_1 < \tau_2$,
	\eqref{eq:SOCWE} is stronger than \eqref{eq:SOCLM}.
	In fact, \eqref{eq:SOCWE} is equivalent to \eqref{eq:SOCLM}
	and $(\nabla F(u))'(1) \ge 0$.
	
	It can be checked that the second order conditions
	obtained via the full representation of \cref{lem:structure}
	are always stronger (or equivalent)
	to the second order conditions via
	the minimal representation (\cite[Corollary~4.4]{LeyfferManns2021}).
	This is also expected
	if we compare \cref{thm:equivalence}
	with the corresponding result
	\cite[Theorem~4.14~(3)]{LeyfferManns2021}.
\end{remark}
We generalize the findings of this example.

\begin{lemma}
	\label{lem:compare_second_order}
	Let $u \in \BV(0,T)$ be feasible for \eqref{eq:prob}
	and we denote by
	$(n,a,t)$ and $(\hat n,\hat a, \hat t)$
	the minimal and the full representation of $u$, respectively.
	We assume that $F$ satisfies the regularity assumptions from \cref{thm:second_derivatives}.
	Further, we define
	the symmetric matrices
	$\F \in \R^{(n-1) \times (n-1)}$
	and $\hat\F \in \R^{(\hat n-1) \times (\hat n-1)}$
	via
	\begin{align*}
		\tau^\top \F \tau
		&:=
		\sum_{j = 1}^{n-1} \mu_j (\nabla F(u))'(t_j) \tau_j^2
		+
		\sum_{j,k = 1}^{n-1} \mu_j \mu_k \nabla^2 F(u)(t_j, t_k) \tau_j \tau_k
		\qquad\forall \tau \in \R^{n-1},
		\\
		\hat\tau^\top \hat\F \hat\tau
		&:=
		\sum_{i = 1}^{\hat n-1} \hat\mu_i (\nabla F(u))'(\hat t_i) \hat\tau_i^2
		+
		\sum_{i,l = 1}^{\hat n-1} \hat\mu_i \hat\mu_l \nabla^2 F(u)(\hat t_i, \hat t_l) \hat\tau_i \hat\tau_l
		\qquad\forall \hat\tau \in \R^{\hat n-1},
	\end{align*}
	where
	$\mu_j = a_{j+1} - a_j$
	and
	$\hat\mu_j = \hat a_{j+1} - \hat a_j$.
	Further, we define the cone
	\begin{equation*}
		\hat\TT :=
		\set[\big]{
			\hat\tau \in \R^{\hat n-1}
			\given
			\forall k \in \set{1,\ldots,\hat n-2}
			:
			\hat t_k = \hat t_{k+1}
			\;\Rightarrow\;
			\hat\tau_k \le \hat\tau_{k+1}
		}
		.
	\end{equation*}
	Then,
	\begin{align*}
		\hat\tau^\top \hat\F \hat\tau \ge 0
		\quad\forall \hat\tau \in \hat\TT
		\qquad\Leftrightarrow\qquad
		\F \succeq 0
		&\;\land\;
		(\nabla F(u))'(t_j) \ge 0 \quad\forall j \in J^+
		\\
		&\;\land\;
		(\nabla F(u))'(t_j) \le 0 \quad\forall j \in J^-
	\end{align*}
	and
	\begin{align*}
		\hat\tau^\top \hat\F \hat\tau > 0
		\quad\forall \hat\tau \in \hat\TT\setminus\set{0}
		\qquad\Leftrightarrow\qquad
		\F \succ 0
		&\;\land\;
		(\nabla F(u))'(t_j) > 0 \quad\forall j \in J^+
		\\
		&\;\land\;
		(\nabla F(u))'(t_j) < 0 \quad\forall j \in J^-
		.
	\end{align*}
\end{lemma}
\begin{proof}
	For $j \in \set{1,\ldots,n}$,
	we set
	$I_j := \set{ i \in \set{1,\ldots,\hat n} \given \hat t_i = t_j}$.
	Note that these sets $I_j$ are a decomposition of $\set{1,\ldots,\hat n-1}$
	and
	$\sum_{i \in I_j} \hat\mu_i = \mu_j$.
	Further, $I_j$ is a singleton,
	if and only if $j \not\in J^+ \cup J^-$.
	Now, let
	$\tau \in \R^{n-1}$
	and $\hat\tau \in \R^{\hat n - 1}$
	be given such that
	\begin{equation}
		\label{eq:rel_tau_hat_tau}
		\tau_j =
		\frac
		{\sum_{i \in I_j} \hat\mu_i \hat\tau_i}
		{\sum_{i \in I_j} \hat\mu_i}
		=
		\frac
		{\sum_{i \in I_j} \hat\mu_i \hat\tau_i}
		{\mu_j}
		\qquad\forall j = 1,\ldots, n-1
		.
	\end{equation}
	Then,
	\begin{align*}
		\hat\tau^\top \hat \F \hat\tau
		&=
		\sum_{j = 1}^{n-1} (\nabla F(u))'(t_j)
		\parens[\Bigg]{ \sum_{i \in I_j} \hat\mu_i \hat\tau_i^2 }
		+
		\sum_{j,k = 1}^{n-1} \nabla^2 F(u)(t_j, t_k)
		\parens[\Bigg]{ \sum_{i \in I_j} \hat\mu_i \hat\tau_i }
		\parens[\Bigg]{ \sum_{l \in I_k} \hat\mu_l \hat\tau_l }
		\\&
		=
		\sum_{j = 1}^{n-1} (\nabla F(u))'(t_j)
		\parens[\Bigg]{ \sum_{i \in I_j} \hat\mu_i \hat\tau_i^2 - \mu_j \tau_j^2}
		+
		\tau^\top \F \tau.
	\end{align*}
	If $I_j$ is a singleton, the last parenthesis vanishes.
	Otherwise,
	\begin{align*}
		\sum_{i \in I_j} \hat\mu_i \hat\tau_i^2
		=
		\mu_j 
		\sum_{i \in I_j} \frac{\hat\mu_i \hat\tau_i^2}{\mu_j}
		\ge
		\mu_j \tau_j^2
		\qquad\forall j \in J^+,
		\\
		\sum_{i \in I_j} \hat\mu_i \hat\tau_i^2
		=
		\mu_j 
		\sum_{i \in I_j} \frac{\abs{\hat\mu_i} \hat\tau_i^2}{\abs{\mu_j}}
		\le
		\mu_j \tau_j^2
		\qquad\forall j \in J^-,
	\end{align*}
	where we used convexity of $s \mapsto s^2$, $\mu_j = \sum_{i\in I_j} \hat\mu_i$
	and that all $\hat\mu_i$, $i \in I_j$ possess the same sign as $\mu_j$.
	Thus,
	\begin{equation}
		\label{eq:rel_F_hat_F}
		\hat\tau^\top \hat\F \hat\tau
		=
		\sum_{j \in J^+ \cup J^-} (\nabla F(u))'(t_j) \sigma_j + \tau^\top \F\tau,
	\end{equation}
	where $\sigma_j = \sum_{i \in I_j} \hat\mu_i \hat\tau_i^2 - \mu_j \tau_j^2$.
	Note that $\pm\sigma_j \ge 0$ for all $j \in J^{\pm}$.

	``$\Rightarrow$'':
	Let $\tau \in \R^{n-1}$ be given.
	We set $\hat\tau_i := \tau_j$ for all $i \in I_j$, $j = 1,\ldots,n-1$.
	Thus, $\hat\tau \in \hat\TT$, \eqref{eq:rel_tau_hat_tau} is satisfied and $\sigma_j = 0$ for all $j \in J^+ \cup J^-$.
	Hence, the positive (semi)-definiteness of $\F$ follows from \eqref{eq:rel_F_hat_F}.

	In order to get the sign conditions of $(\nabla F(u))'(t_j)$ for $j \in J^+ \cup J^-$,
	it is enough to realize that we can choose $\hat\tau \in \hat\TT$
	such that the corresponding $\tau$ and $\sigma$
	satisfy $\tau = 0$, $\sigma_j = \pm 1$ and $\sigma_{\tilde j} = 0$
	for $\tilde j \in (J^+ \cup J^-) \setminus \set{j}$.

	``$\Leftarrow$'':
	For a given $\hat\tau \in \hat\TT$,
	let $\tau$ according to
	\eqref{eq:rel_tau_hat_tau}
	be given.
	Then, $\hat\tau^\top \hat\F \hat\tau \ge 0$ follows from \eqref{eq:rel_F_hat_F}.

	It remains to prove the positive definiteness under the stronger conditions
	on $\F$ and $(\nabla F(u))'(t_j)$.
	One can check that $\hat\tau \ne 0$
	implies
	$\tau \ne 0$ or $\sigma_j \ne 0$ for some $j \in J^+ \cup J^-$.
	Thus, we get
	$\hat\tau^\top \hat\F \hat\tau > 0$ from \eqref{eq:rel_F_hat_F}.
\end{proof}
Note that the conditions involving the cone $\hat\TT$
are difficult to verify
since they involve
positive (semi)-definiteness of a matrix
over a cone
and this is, in general, difficult to check.
In contrast, the equivalent conditions appearing on the right-hand sides
are straightforward to verify.

By combining the above results,
we obtain the main result of this section.
\begin{theorem}
	\label{thm:no_gap_SOC}
	Let $u \in \BV(0,T)$ be feasible for \eqref{eq:prob}
	and we denote by $(n,a,t)$
	the minimal representation for $u$.
	We assume that $F \colon L^1(0,T) \to \R$
	is twice Fréchet differentiable
	with $\nabla F(u) \in C^1([0,T])$
	and $\nabla^2 F(u) \in C([0,T]^2)$.
	We define
	$\mu_j := a_{j+1} - a_j$ for $j = 1,\ldots,n-1$.
	If $u$ is a local minimizer of \eqref{eq:prob} in $L^1(0,T)$, then the system
	\begin{subequations}
		\label{eq:SONC}
		\begin{align}
			\label{eq:SONC_1}
			\nabla F(u)(t_j) &= 0 \qquad\forall j = 1,\ldots, n-1, \\
			\label{eq:SONC_2}
			(\nabla F(u))'(t_j) &\ge 0 \qquad\forall j \in J^+, \\
			\label{eq:SONC_3}
			(\nabla F(u))'(t_j) &\le 0 \qquad\forall j \in J^-, \\
			\label{eq:SONC_4}
			\sum_{j = 1}^{n-1} \mu_j (\nabla F(u))'(t_j) \tau_j^2
			+
			\sum_{j,k = 1}^{n-1} \mu_j \mu_k \nabla^2 F(u)(t_j, t_k) \tau_j \tau_k
			&\ge0
			\qquad\forall \tau \in \R^{n-1}
		\end{align}
	\end{subequations}
	is satisfied.
	Moreover,
	$u$ is a local minimizer of \eqref{eq:prob} satisfying
	a quadratic growth condition in $L^1(0,T)$
	if and only if
	\begin{subequations}
		\label{eq:SOEC}
		\begin{align}
			\label{eq:SOEC_1}
			\nabla F(u)(t_j) &= 0 \qquad\forall j = 1,\ldots, n-1, \\
			\label{eq:SOEC_2}
			(\nabla F(u))'(t_j) &> 0 \qquad\forall j \in J^+, \\
			\label{eq:SOEC_3}
			(\nabla F(u))'(t_j) &< 0 \qquad\forall j \in J^-, \\
			\label{eq:SOEC_4}
			\sum_{j = 1}^{n-1} \mu_j (\nabla F(u))'(t_j) \tau_j^2
			+
			\sum_{j,k = 1}^{n-1} \mu_j \mu_k \nabla^2 F(u)(t_j, t_k) \tau_j \tau_k
			&> 0
			\qquad\forall \tau \in \R^{n-1} \setminus\set{0}.
		\end{align}
	\end{subequations}
\end{theorem}
Note that \eqref{eq:SONC_4}, \eqref{eq:SOEC_4}
describe the positive (semi)-definiteness
of the matrix
\begin{equation*}
	\diag\parens*{\bracks*{
			\mu_j (\nabla F(u))'(t_j)
		}_{j = 1,\ldots,n-1}
	}
	+
	\bracks*{
		\mu_j \mu_k \nabla^2 F(u)(t_j, t_k)
	}_{j,k = 1,\ldots, n-1}.
\end{equation*}
Furthermore, we mention that
\eqref{eq:SONC} and \eqref{eq:SOEC}
can be easily checked.
Bear in mind that these conditions
use the data from the minimal representation of $u$,
but were derived using the full representation of $u$.
Finally, we mention that the gap between the necessary
and the sufficient conditions
is as small as possible
and, moreover,
we are able to characterize local quadratic growth in $L^1(0,T)$.

\begin{remark}\leavevmode
	\label{rem:second_derivatives}
	\begin{enumerate}
		\item
			A comparable second-order optimality condition
			(for bang-bang problems)
			in the multi-dimensional case
			was given in
			\cite[Theorem~6.12]{ChristofWachsmuth2017:1}.
			Therein, the term $\abs{\nabla \varphi}$
			corresponds to $(\nabla F(v^{t,a}))'$ above
			(since the adjoint state $\varphi$ represents the derivative
			of the objective w.r.t.\ the control at the point of interest).
		\item
			The results of \cref{thm:second_derivatives,thm:no_gap_SOC} can be utilized
			to set up a Newton method for the solution of \eqref{eq:ST}.
		\item
			The second-order terms in the \cref{thm:second_derivatives,thm:no_gap_SOC}
			give rise to the following observations:
			\begin{itemize}
				\item
					The convexity of $F$ is not enough to guarantee
					that first order stationary points are (locally) optimal.
					Indeed, the convexity of $F$ has no influence
					on the signs of $(\nabla F(u))'(t_j)$.
				\item
					Similarly, optimality of $v^{t,a}$ alone does
					not give a sign of $(\nabla F(v^{t,a}))'(t_j)$
					for $j \not\in J^+ \cup J^-$,
					due to the coupling in \eqref{eq:SONC_4}.
			\end{itemize}
	\end{enumerate}
\end{remark}

\subsection{Non-local optimality conditions}
\label{subsec:non-local_optimality}

In \cref{thm:no_gap_SOC},
we were able
to give second-order optimality conditions
with minimal gap.
This delivers a good understanding
of the local optimality for the problem \eqref{eq:prob}.

In this section, we provide two examples of a non-local optimality condition.
The first result shows that fast back-and-forth switches
can be non-optimal in certain situations.
\begin{theorem}
	\label{thm:non_optimal_switch}
	Suppose that $F \colon L^1(0,T) \to \R$ is Fréchet differentiable with
	Lipschitz continuous derivative $\nabla F \colon L^1(0,T) \to L^\infty(0,T)$
	and Lipschitz constant $L \ge 0$.
	Further, let $u \in \BV(0,T)$ be feasible for \eqref{eq:prob}
	and let $j \in \set{1,\ldots,d}$ and $0 < t_1 < t_2 < t_3 < t_4 < T$ be given,
	such that $u = \nu_j$ on $(t_2, t_3)$ and
	$u < \nu_j$ on $(t_1,t_2) \cup (t_3, t_4)$ hold.
	If
	\begin{equation}
		\label{eq:non_optimal_switch}
		\mathopen{}- 2 \beta
		- \int_{t_2}^{t_3} \nabla F(u)(s) \, \ds
		+
		\frac{L}{2} (\nu_j - \nu_{j-1}) (t_3 - t_2)^2
		<
		0
	\end{equation}
	then $v = u + (\nu_{j-1} - \nu_j) \chi_{(t_2, t_3)}$
	satisfies
	\begin{equation*}
		F(v) + \beta \TV(v)
		<
		F(u) + \beta \TV(u).
	\end{equation*}
\end{theorem}
\begin{proof}
	The Lipschitz continuity of $\nabla F$ ensures
	\begin{equation*}
		\abs{
			F(v) - F(u) - \nabla F(u)(v - u)
		}
		\le \frac{L}{2} \norm{v - u}_{L^1(0,T)}^2
		=
		\frac{L}{2} (\nu_j - \nu_{j-1})^2 (t_3 - t_2)^2
		.
	\end{equation*}
	Together with
	\begin{equation*}
		\nabla F(u)(v - u)
		=
		(\nu_{j-1} - \nu_j) \int_{t_2}^{t_3} \nabla F(u)(s) \, \ds
	\end{equation*}
	and
	\begin{equation*}
		\TV(v) = \TV(u) - 2 (\nu_j - \nu_{j-1})
		,
	\end{equation*}
	this establishes the claim.
\end{proof}
Note that \cref{thm:non_optimal_switch}
is concerned with the situation of $u$ switching upwards on $(t_2, t_3)$.
A similar argument can be used in case of a downward switch with
$u > \nu_j$ on $(t_1,t_2) \cup (t_3, t_4)$.

We mention that \eqref{eq:non_optimal_switch}
is always satisfied if $t_3 - t_2$ is small enough.
Indeed,
if $\abs{\nabla F(u)} \le C$ holds on $(0,T)$,
then
\begin{equation*}
	t_3 - t_2
	<
	\frac{-C + \sqrt{C^2 + 4 \beta L (\nu_j - \nu_{j-1})}}{L (\nu_j - \nu_{j-1})}
\end{equation*}
implies \eqref{eq:non_optimal_switch}.

Finally, we comment that $u$ can still be locally optimal in the situation of
\cref{thm:non_optimal_switch}.
To see this, consider that $\norm{u - v}_{L^1(0,T)} = (\nu_j - \nu_{j-1}) (t_3 - t_2)$
and the radius of optimality of $u$ could be smaller than this constant.

The next result is concerned with the introduction of an additional switch.
\begin{theorem}
	\label{thm:non_optimal_non_switch}
	Suppose that $F \colon L^1(0,T) \to \R$ is Fréchet differentiable with
	Lipschitz continuous derivative $\nabla F \colon L^1(0,T) \to L^\infty(0,T)$
	with constant $L \ge 0$.
	Further, let $u \in \BV(0,T)$ be feasible for \eqref{eq:prob}
	and let $j \in \set{1,\ldots,d}$ and $0 < t_1 < t_4 < T$ be given,
	such that $u = \nu_j$ on $(t_1, t_4)$.
	Suppose that
	\begin{equation}
		\label{eq:non_optimal_non_switch}
		2 \beta \abs{\nu_k - \nu_j}
		+ (\nu_k - \nu_j) \int_{t_2}^{t_3} \nabla F(u)(s) \, \ds
		+
		\frac{L}{2} (\nu_k - \nu_j)^2 (t_3 - t_2)^2
		<
		0
	\end{equation}
	is satisfied,
	where
	$t_1 < t_2 < t_3 < t_4$
	and $k \in \set{1, \ldots, d} \setminus \set{j}$.
	Then $v = u + (\nu_k - \nu_j) \chi_{(t_2, t_3)}$
	satisfies
	\begin{equation*}
		F(v) + \beta \TV(v)
		<
		F(u) + \beta \TV(u).
	\end{equation*}
\end{theorem}
\begin{proof}
	This follows from similar arguments as in the proof of \cref{thm:non_optimal_switch},
	but now we have
	\begin{align*}
		\abs{
			F(v) - F(u) - \nabla F(u)(v - u)
		}
		&\le \frac{L}{2} \norm{v - u}_{L^1(0,T)}^2
		=
		\frac{L}{2} (\nu_k - \nu_j)^2 (t_3 - t_2)^2
		,\\
		\nabla F(u)(v - u)
		&=
		(\nu_k - \nu_j) \int_{t_2}^{t_3} \nabla F(u)(s) \, \ds
		,\\
		\TV(v) &= \TV(u) + 2 \abs{\nu_k - \nu_j}
		.
	\end{align*}
\end{proof}
This result shows that it might be worthwhile to
have jumps to bigger/smaller values
when
$\nabla F(u)$
is negative/positive
on intervals where $u$ is constant.
In contrast to \cref{thm:non_optimal_non_switch},
the region $(t_2,t_3)$ on which $u$ will be modified
cannot be too small, otherwise the first term in \eqref{eq:non_optimal_non_switch} dominates.

\section{Proximal-gradient method}
\label{sect:Proximal_Gradient}
In this section, we propose a proximal-gradient method to compute locally optimal points of \eqref{eq:prob}.
Originally, the method was proposed for non-differentiable convex optimization problems,
but
contributions like \cite{Wachsmuth2019:3} motivate the application
to non-convex problems,
also in infinite dimensions.

\subsection{Theoretical results}
\label{subsec:prox_theory}
Since the proximal-gradient method applies to problems in Hilbert spaces,
we will discuss \eqref{eq:prob}
in the space $L^2(0,T)$.
Note that the admissible set $\Uad$ is already a subset of $L^2(0,T)$.
We start by reformulating \eqref{eq:prob}
as
 \[
	\min_{u\in L^2(0,T)}  F(u) +\beta\TV(u)+\delta_{\Uad}(u)
 \]
where the indicator function $\delta_{\Uad}\colon L^1(0,T)\to \set{0,\infty}$ is defined by
\[
	\delta_{\Uad} = \begin{cases}
		0,\quad&\text{if } u\in \Uad,\\
		\infty,\quad&\text{otherwise}.
	\end{cases}
\]
Now, the first addend in the objective $F$ is smooth,
whereas the second part
$G\colon L^2(0,T)\to \R\cup \{\infty\}$, given by
\[
	G(u) := \beta\TV(u)+\delta_{\Uad}(u),
\]
is non-smooth and non-convex.
As in \cite[Algorithm~3.21]{Wachsmuth2019:3},
we use
the decrease condition
\begin{equation}
	\label{eq:decrease_condition}
	\eta\norm{u_{k+1}-u_k}^2_{L^2(0,T)}\le F(u_k) + \beta\TV(u_k) - (F(u_{k+1}) + \beta\TV(u_{k+1})),
\end{equation}
with some parameter $\eta>0$
in each step of the proximal-gradient method,
see
\cref{Alg:prox-grad}.
\begin{algorithm2e}[tb]
	\SetAlgoLined
	\KwData{$F,G\colon L^2(0,T)\to \R$, where $F$ is Gateaux-differentiable, $u_0\in \Uad$, $\eta>0$}
	\KwResult{Sequence $\set{u_k}_{k\in\N}\subset L^2(0,T)$}
	Choose $\tau_k>0$ such that a solution $u_{k+1}$ of
	\begin{equation}
		\label{eq:prox_grad_iter}
		\min_{u\in L^2(0,T)}  F(u_k) + \nabla F(u_k)(u-u_k) + \frac{\tau_k}{2} \norm{u-u_k}^2_{L^2(0,T)} + G(u)
	\end{equation}
	satisfies \eqref{eq:decrease_condition}.
	
	Set $k\gets k+1$ and go to step 1.
	\caption{Proximal-Gradient Algorithm}
	\label{Alg:prox-grad}
\end{algorithm2e}
The existence of solutions $u_{k+1}$ of problem \eqref{eq:prox_grad_iter}
can be guaranteed similar to the discussion after \cref{thm:props}.
However, since $G$ fails to be convex, there might be multiple solutions.
The next result gives some basic properties of sequences generated by \cref{Alg:prox-grad}.
\begin{theorem}\label{thm:prox}
	Let $\seq{u_k}_{k\in\N}$ be a sequence generated by \cref{Alg:prox-grad}.
	Moreover, let $\nabla F$ be Lipschitz continuous from $L^2(0,T)$ to $L^2(0,T)$ with modulus $L$.
	Then, the following is true:
	\begin{enumerate}
		\item\label{thm:prox:1} The sequences $\seq{u_k}_{k\in\N}$ and $\seq{\nabla F(u_k)}_{k\in\N}$ are bounded in $L^2(0,T)$.
		\item\label{thm:prox:2} The sequence $\seq{F(u_k)+G(u_k)}_{k\in\N}$ is decreasing and converges.
		\item\label{thm:prox:3} $\norm{u_{k+1}-u_k}_{L^2(0,T)}\rightarrow 0$.
		\item\label{thm:prox:4} $\seq{u_k}_{k\in\N}$ converges weak-$\star$ in $\BV(0,T)$ towards some $\bar u \in \Uad$.
	\end{enumerate}
\end{theorem}
\begin{proof}
	We will adapt the proof of \cite[Theorem 3.22]{Wachsmuth2019:3} for our situation. Since \eqref{eq:decrease_condition} can be written as
	\[F(u_{k+1}) + G(u_{k+1})\leq F(u_k) + G(u_k) - \eta\norm{u_{k+1}-u_k}_{L^2(0,T)}^2,\]
	and $F$, $G$ are bounded from below, \ref{thm:prox:2} follows.
	This implies that $\seq{G(u_k)}_{k\in\N}$ is also bounded.
	Furthermore, as we have $G(u) = \infty$ for $u\notin\Uad$, $u_k\in\Uad$ holds for all $k\in\N$. Thus, 
	\begin{equation*}
		\norm{u_k}_{L^2(0,T)}^2 \le T \max\set{\abs{\nu_1}, \abs{\nu_d}}^2.
	\end{equation*}
	Moreover, using the Lipschitz continuity
	of $\nabla F$, this implies the boundedness of $\seq{\nabla
	F(u_k)}$, which completes the proof of \ref{thm:prox:1}.

	By taking the sum of \eqref{eq:decrease_condition} over $k=1,\dots,n$ for $n\in\N$ leads to 
	\begin{equation*}
		\begin{aligned}
			F(u_{n+1}) + G(u_{n+1}) + \eta\sum_{k=1}^n \norm{u_{k+1}-u_k}^2_{L^2(0,T)}&\leq F(u_1) + G(u_1).
		\end{aligned}
	\end{equation*}
	With $n\to\infty$, we see that
	\[\lim_{n\to\infty}\parens*{ F(u_{n+1}) + G(u_{n+1})} + \eta\sum_{k=1}^\infty \norm{u_{k+1}-u_k}^2_{L^2(0,T)}\leq F(u_1) + G(u_1) < \infty,\]
	which implies that the series $\sum_{k=1}^\infty \norm{u_{k+1}-u_k}^2_{L^2(0,T)}$ converges.
	Thus, \ref{thm:prox:3} follows. 

	To show \ref{thm:prox:4}, we note that $\abs{u_{k+1}-u_k}$ does not take values 
	in $(0,1)$ for all $k\in\N$.
	This leads to the inequality 
	$\abs{u_{k+1}-u_k}^2\geq\abs{u_{k+1}-u_k}$ and hence 
	$\norm{u_{k+1}-u_k}_{L^1(0,T)}\leq \norm{u_{k+1}-u_k}_{L^2(0,T)}^2$. Now, 
	since $\sum_{k=1}^\infty \norm{u_{k+1}-u_k}^2_{L^2(0,T)}$ converges, 
	we get the convergence of the series
	$\sum_{k=1}^\infty \norm{u_{k+1}-u_k}_{L^1(0,T)}$, which leads to the 
	strong convergence of $\seq{u_k}_{k\in\N}$ in $L^1(0,T)$.
	Since $F$ is bounded from below,
	the sequence $\seq{u_k}$ is bounded in $\BV(0,T)$.
	This shows
	$u_k \weaklystar \bar u$ in $BV(0,T)$,
	see \itemref{thm:props:1}.
	Finally, since $\Uad$ is closed in $L^1(0,T)$,
	$\bar u \in \Uad$ follows.
\end{proof}
Note that \cite[Theorem 3.13]{Wachsmuth2019:3} states the validity of 
\[F(u_{k+1})+G(u_{k+1})\leq F(u_k) + G(u_k)-\frac{\tau_k-L}{2}\norm{u_{k+1}-u_k}_{L^2(0,T)}^2,\]
where $u_{k+1}$ is the solution of \eqref{eq:prox_grad_iter}.
Hence the choice $\tau_k \ge 2\eta + L$
implies that
the decrease condition \eqref{eq:decrease_condition} is satisfied.
Nevertheless, for a fast 
convergence of the algorithm, it is desired to choose the inverse step length $\tau_k$ as small as possible.
This can be
realized by testing the values $\tau^0\theta^{-i}$ for $i=0,1,2,\dots$,
$\tau^0>0$ and $\theta\in(0,1)$ until the decrease condition is achieved. If
$\tau^0$ is already sufficient, it is reasonable to test smaller values
$\tau^0\theta^i$ for $i=1,2,\dots$ until \eqref{eq:decrease_condition} is no longer
valid.

\begin{theorem}
	\label{thm:prox_opt}
	Let $\seq{u_k}_{k\in\N}$ be a sequence generated by \cref{Alg:prox-grad}.
	Further, let $\nabla F$ be Lipschitz continuous from $L^2(0,T)$ to $L^2(0,T)$ with modulus $L$.
	Then,
	the weak-$\star$ limit $\bar u$
	of the sequence $\seq{u_k}_{k\in\N}$ in $\BV(0,T)$
	solves
	\begin{equation}
		\label{eq:prox_opt}
		\min_{u\in \Uad}  F(\bar u) + \nabla F(\bar u)(u-\bar u) + \frac{\bar\tau}{2} \norm{u-\bar u}^2_{L^2(0,T)} + \beta \TV(u)
	\end{equation}
	for every accumulation point $\bar \tau$ of $\seq{\tau_k}$.
\end{theorem}
\begin{proof}
	Since $u_{k + 1}$ solves \eqref{eq:prox_grad_iter},
	we have
	\begin{align*}
		&
		\nabla F(u_{k})(u_{k+1}-u_{k}) + \frac{\tau_{k}}{2} \norm{u_{k+1}-u_{k}}_{L^2(0,T)}^2 + \beta \TV(u_{k+1})
		\\&\qquad
		\le
		\nabla F(u_{k})(v        -u_{k}) + \frac{\tau_{k}}{2} \norm{v        -u_{k}}_{L^2(0,T)}^2 + \beta \TV(v        )
	\end{align*}
	for all $v \in L^2(0,T) \cap \Uad = \Uad$.
	Suppose that the subsequence $\seq{\tau_{k_l}}$ converges towards $\bar \tau$.
	The above inequality yields
	\begin{align*}
		\beta \TV(\bar u)
		&
		\le
		\liminf_{l \to \infty}
		\parens*{
			\nabla F(u_{k_l})(u_{k_l+1}-u_{k_l}) + \frac{\tau_{k_l}}{2} \norm{u_{k_l+1}-u_{k_l}}_{L^2(0,T)}^2 + \beta \TV(u_{k_l+1})
		}
		\\&
		\le
		\lim_{l \to \infty}
		\parens*{
			\nabla F(u_{k_l})(v        -u_{k_l}) + \frac{\tau_{k_l}}{2} \norm{v        -u_{k_l}}_{L^2(0,T)}^2 + \beta \TV(v        )
		}
		\\&
		=
		\nabla F(\bar u)(v-\bar u) + \frac{\bar\tau}{2} \norm{v -\bar u}_{L^2(0,T)}^2 + \beta \TV(v)
		.
	\end{align*}
	Since $v \in \Uad$ was arbitrary, this shows the claim.
\end{proof}
Next, we are going to investigate optimality conditions of \eqref{eq:prox_opt}.
Note that it is not possible to utilize the theory of
\cref{sec:optimality_conditions},
since $u \mapsto \frac{\bar\tau}{2}\norm{u - \bar u}_{L^2(0,T)}^2$
is not Fréchet differentiable in $L^1(0,T)$.
The following lemma shows
that the optimality conditions of \eqref{eq:prox_opt}
are weaker than the first order conditions from \cref{thm:no_gap_SOC}.
\begin{lemma}
	\label{lem:prox_opt_con}
	Let $\bar u \in \BV(0,T) \cap \Uad$ and $\bar\tau \ge 0$ be given
	such that $\bar u$ is a solution of \eqref{eq:prox_opt}.
	Further, suppose that $\nabla F(\bar u) \in C([0,T])$.
	Then, for each switching time $t \in (0,T)$,
	we have
	\begin{subequations}
		\label{eq:prox_opt_con}
		\begin{align}
			\hat a_i &< \hat a_{j+1}
			\qquad\Rightarrow\qquad
			-\frac{\bar\tau}{2} \abs{\hat a_{i+1} - \hat a_i}
			\le
			(\nabla F(\bar u))(t)
			\le
			\frac{\bar\tau}{2} \abs{\hat a_{j+1} - \hat a_j}
			,
			\\
			\hat a_i &> \hat a_{j+1}
			\qquad\Rightarrow\qquad
			-\frac{\bar\tau}{2} \abs{\hat a_{j+1} - \hat a_j}
			\le
			(\nabla F(\bar u))(t)
			\le
			\frac{\bar\tau}{2} \abs{\hat a_{i+1} - \hat a_i}
			,
		\end{align}
	\end{subequations}
	in which
	we use the data $(\hat t, \hat a)$ from the full representation,
	$i$ is the smallest index with $t = \hat t_i$
	and
	$j$ is the largest index with $t = \hat t_j$.
\end{lemma}
In the case that $\nu_{i+1} - \nu_i=1$ for all $i = 1, \ldots, d-1$,
\eqref{eq:prox_opt_con} is equivalent to
$ \abs{(\nabla F(\bar u))(t)} \le \bar\tau/2 $.
\begin{proof}
	For an arbitrary $\varepsilon \in (0, \hat t_i - \hat t_{i-1})$, we consider the perturbed function
	\begin{equation*}
		v_\varepsilon := \bar u + (\hat a_{i+1} - \hat a_i) \chi_{(t-\varepsilon, t)},
	\end{equation*}
	i.e., we change the value of $\bar u$ on $(t - \varepsilon, t)$
	from $\hat a_i$ to $\hat a_{i+1}$.
	Thus, $\TV(v_\varepsilon) = \TV(\bar u)$
	and the optimality of $\bar u$ for \eqref{eq:prox_opt} gives
	\begin{align*}
		0
		&\le
		\nabla F(\bar u)(v_\varepsilon - \bar u) + \frac{\bar\tau}{2} \norm{v_\varepsilon - \bar u}_{L^2(0,T)}^2
		\\&=
		(\hat a_{i+1} - \hat a_i)
		\int_{t - \varepsilon}^{t}
		(\nabla F(\bar u))(s) \, \d s
		+
		\frac{\bar\tau}{2} \abs{\hat a_{i+1} - \hat a_i}^2 \varepsilon
		.
	\end{align*}
	Dividing by $\varepsilon > 0$ and passing to the limit $\varepsilon \searrow 0$
	yields
	\begin{equation*}
		0
		\le
		(\hat a_{i+1} - \hat a_i)
		(\nabla F(\bar u))(t)
		+
		\frac{\bar\tau}{2} \abs{\hat a_{i+1} - \hat a_i}^2
		.
	\end{equation*}
	Similarly, we can use the perturbation
	\begin{equation*}
		v_\varepsilon := \bar u + (\hat a_{j} - \hat a_{j+1}) \chi_{(\hat t_j, \hat t_j+\varepsilon)}
	\end{equation*}
	and this leads to
	\begin{equation*}
		0
		\le
		(\hat a_{j} - \hat a_{j+1})
		(\nabla F(\bar u))(t)
		+
		\frac{\bar\tau}{2} \abs{\hat a_{j} - \hat a_{j+1}}^2
		.
	\end{equation*}
	Using the observation that the signs of $\hat a_{j+1} - \hat a_i$, $\hat a_{i+1} - \hat a_i$
	and $\hat a_{j+1} - \hat a_j$ coincide,
	see \itemref{lem:structure:4},
	we arrive at \eqref{eq:prox_opt_con}.
\end{proof}
Note that the condition \eqref{eq:prox_opt_con}
is weaker than the first-order condition \eqref{eq:SONC_1}
in case $\bar\tau > 0$.
A similar observation has been made in \cite[Theorem~3.18]{Wachsmuth2019:3}.

\subsection{Fast solution of discrete subproblems}
\label{sect:prox_subproblems}
The main work of \cref{Alg:prox-grad}
consists in the solution of
the subproblems
\eqref{eq:prox_grad_iter},
which can be equivalently written as
\begin{equation}
	\label{eq:subproblem1}
	\min_{u\in \Uad} F(u_k) + \nabla F(u_k)(u-u_k) + \frac{\tau_k}{2}\norm{u-u_k}^2_{L^2(0,T)} + \beta\TV(u)
	.
\end{equation}
On a first glance,
these subproblems seem to be very
delicate,
since we have the integer constraints,
some nonlinearity and
the coupling in time due to the $\TV$-norm.
However, we will see
that it is possible to solve (the discretizations of) these problems
very efficiently.

First, we want to restate \eqref{eq:subproblem1}.
We define the gradient step
\[v_k:=u_k - \frac{1}{\tau_k}\nabla F(u_k)\in L^2(0,T).\]
We can use
\begin{equation*}
\frac{\tau_k}{2}\norm{u-v_k}^2_{L^2(0,T)}=\frac{\tau_k}{2}\norm{u-u_k}_{L^2(0,T)}^2
		+\nabla F(u_k)(u-u_k)
		+\frac{1}{2\tau_k}\norm{\nabla F(u_k)}_{L^2(0,T)}^2,
\end{equation*}
to rewrite the objective of \eqref{eq:subproblem1}.
By further
omitting the constant terms and by dropping the index $k$ of $v_k$ and $\tau_k$,
\eqref{eq:subproblem1} can be rephrased as
\begin{equation}
	\label{eq:subproblem2}
	\min_{u\in \Uad} \frac{\tau}{2}\norm{u-v}^2_{L^2(0,T)} + \beta\TV(u).
\end{equation}
Note that the solution of \eqref{eq:subproblem2}
corresponds to the computation of the proximal point mapping
of the non-convex functional
$G = \beta \TV + \delta_{\Uad}$.

In order to discretize
\eqref{eq:subproblem2},
we partition $[0,T]$ via the grid $0=t_0<t_1<\dots<t_n = T$.
For simplicity of the presentation,
we assume that we have
an equidistant mesh size $\Delta t := \frac{T}{n}$,
but the following
can be adapted easily to non-equidistant mesh sizes.

In accordance with this mesh, we discretize the function $u$
as a piecewise constant function,
i.e.,
$u = \sum_{j = 1}^n u^j \chi_{(t_{j-1}, t_j)}$,
for $u^j \in \set{\nu_1, \ldots, \nu_d}$, $j = 1,\ldots, n$.
For the discretization of $v$, we choose the mean values
$v^j = (\Delta t)^{-1} \int_{t_{j-1}}^{t_j} v \, \dt$.
Thus, a discretization of \eqref{eq:subproblem2} is given by
\begin{equation}
	\label{eq:subproblem3}
	\min_{u^1,\dots,u^{n}\in\set{\nu_1,\dots,\nu_d}}
	\frac{\tau \Delta t}{2}\sum_{j=1}^{n} (u^j - v^j)^2
	+ \beta \sum_{j=1}^{n-1} \abs{u^{j+1}-u^{j}}
\end{equation}
or, equivalently,
\begin{equation}
	\label{eq:subproblem4}
	\min_{\kappa_1,\dots,\kappa_{n}\in\set{1,\dots,d}}
	\frac{\tau \Delta t}{2}\sum_{j=1}^{n} (\nu_{\kappa_j} - v^j)^2
	+ \beta \sum_{j=1}^{n-1} \abs{\nu_{\kappa_{j+1}}-\nu_{\kappa_j}}.
\end{equation}
Now, we want to employ the Bellman principle on problem \eqref{eq:subproblem4},
stating that independent from the initial decision, the remaining decisions of
an optimal solution have to constitute an optimal policy with regard to the
state resulting from the first decision.
In this sense, we define a value function
(represented by the matrix $\Phi \in \R^{d \times n}$)
giving the optimal value of
\eqref{eq:subproblem4} restricted to an interval $(t_{\iota-1},T)$ given the
choice $u^\iota = \nu_{\kappa_l}$ at $t_{\iota-1}$.
That is, we define
\begin{equation}
	\label{eq:value_func1}
	\Phi_{l, \iota}
	:=
	\min
	\set*{
		\frac{\tau \Delta t}{2}
		\sum_{j=\iota}^{n} (\nu_{\kappa_j} - v^j)^2
		+ \beta \sum_{j=\iota}^{n-1} \abs{\nu_{\kappa_{j+1}}-\nu_{\kappa_{j}}}
		\given
		\begin{aligned}
			& \kappa_\iota = l,\\
			& \kappa_{\iota+1},\dots,\kappa_n\in\{1,\dots,d\}
		\end{aligned}
	}
\end{equation}
for all $l = 1,\ldots, d$, $\iota = 1,\ldots, n$.
It is easy to see that
this gives
\begin{equation}
	\label{eq:value_func2}
	\forall l =1,\ldots, d:
	\qquad
	\Phi_{l,n} = 
	\frac{\tau \Delta t}{2} (\nu_{l} - v^j)^2,
\end{equation}
which is a terminal value for the value function.
In order to compute $\Phi_{l, \iota}$ for $\iota < n$,
we have to minimize
\begin{equation*}
	\bracks*{
		\frac{\tau \Delta t}{2}
		(\nu_{l} - v^\iota)^2
		+ \beta \abs{\nu_{\kappa_{\iota+1}}-\nu_{l}}
	}
	+
	\bracks*{
		\frac{\tau \Delta t}{2}
		\sum_{j=\iota+1}^{n} (\nu_{\kappa_j} - v^j)^2
		+ \beta \sum_{j=\iota+1}^{n-1} \abs{\nu_{\kappa_{j+1}}-\nu_{\kappa_{j}}}
	}
\end{equation*}
w.r.t.\ $\kappa_{\iota+1}, \ldots, \kappa_n \in \set{1,\ldots,d}$.
The first bracket is independent of $\kappa_{\iota+2}, \ldots, \kappa_n$,
hence, these values minimize the second bracket
and the corresponding minimal value is $\Phi_{\kappa_{\iota+1},\iota+1}$.
Thus, for all $1 \le l \le d$ and $1 \le \iota < n$, \eqref{eq:value_func1} can be rephrased as
\begin{equation}
	\label{eq:value_func3}
	\Phi_{l,\iota}
	=
	\min\set*{
	\frac{\tau \Delta t}{2}
	(\nu_{l} - v^\iota)^2
	+ \beta \abs*{\nu_{\kappa_{\iota+1}}-\nu_{l}}
	+
	\Phi_{\kappa_{\iota+1},\iota+1}
	\given
	\kappa_{\iota+1} \in \set{1,\ldots,d}
	}.
\end{equation}
Finally,
the solution of \eqref{eq:subproblem4} can be found by calculating $\Phi_{l,1}$ for every $l\in\{1,\dots,d\}$ and comparing these values.
As motivated before, this can be achieved by computing $\Phi_{l,\iota}$ for $\iota= n,\dots,1$ and every $l\in\{1,\dots,d\}$ using \eqref{eq:value_func2} in the first step (which, in our case, is the last time step) and \eqref{eq:value_func3} for the following steps. 
The corresponding minimizer $\kappa_{\iota+1}$ has to be saved for every $\iota=n-1,\dots,1$ in order to reconstruct the solution when the best initial choice $l\in\{1,\dots,d\}$ minimizing $\Phi_{l,1}$ has been found. Therefore, we save these values in a matrix $U\in \R^{d\times n-1}$ defined by 
\begin{equation*}
	U_{l,\iota} :=  \argmin\set*{
		\frac{\tau \Delta t}{2}
		(\nu_{l} - v^\iota)^2
		+ \beta \abs*{\nu_{\kappa_{\iota+1}}-\nu_{l}}
		+
		\Phi_{\kappa_{\iota+1},\iota+1}
		\given
		\kappa_{\iota+1} \in \set{1,\ldots,d}
	}
\end{equation*}
for every time step $\iota\in\{1,\dots,n-1\}$.

Now, $u$ can be calculated by setting $\kappa_1:=\argmin\set{\Phi_{l,1}\given l\in\{1,\dots,d\}}$, $\kappa_\iota:= U_{\kappa_{\iota-1},\iota-1}$ for $\iota\in\{2,\dots,n\}$ and $u^\iota = \nu_{\kappa_\iota}$ for $\iota\in\{1,\dots,n\}$.
\begin{remark}
	\label{rem:prox-subproblem}
	In an implementation, only a $d\times 2$ matrix $\Phi$ is needed since we can
	overwrite the old target values in a step $\iota+1$ with the new ones of step
	$\iota$ after $\Phi_{l, \iota}$ has been computed for every $l\in\set{1,\dots,d}$.

	By testing $d$ target values for $d$ possible settings of $l$ and repeating
	this for all $n-1$ time steps, the emerging algorithm has a runtime of
	$\OO(d^2n)$.
\end{remark}

\section{Trust-region algorithm and efficient computation of corresponding subproblems}
\label{sec:TR}

Similar to \cite[Sect. 3.1]{LeyfferManns2021}, locally optimal points of 
\eqref{eq:prob} can be calculated using a trust-region algorithm where the 
objective is partially linearized around a given feasible point. When employing such an algorithm, one has to solve subproblems of the form
\begin{equation}
	\label{eq:TR_subprob}
	\tag{TR}
	\begin{aligned}
		\text{Minimize} \quad & (g,u - v)_{L^2(0,T)} + \beta \TV(u) - \beta\TV(v) \\
		\text{such that} \quad & \norm{u-v}_{L^1(0,T)}\leq \Delta^k,\quad u \in  \Uad
	\end{aligned}
\end{equation}
with a given function $v\in\Uad$ and $g=\nabla F(v)$. In \cite{LeyfferManns2021}, this was done by constructing a mixed-integer linear program. For a fine discretization, such an approach may lead to long computing times, which is why we are interested in applying the Bellman principle in a similar manner as in \cref{sect:Proximal_Gradient} to efficiently compute discrete solutions of \eqref{eq:TR_subprob}.

Therefore, consider the same discretization of $[0,T]$ as in \cref{sect:prox_subproblems} with an equidistant mesh size $\Delta t = \tfrac Tn$ and $v = \sum_{j=1}^n v^j\chi_{(t_{j-1},t_j)}$, $g^j := (\Delta t)^{-1}\int_{t_{j-1}}^{t_j} g\,\d t$, $j=1,\dots,n$. We rephrase the problem by omitting terms in the objective independent of $u$, obtaining the formulation
\begin{equation}
	\label{eq:TR_subprob2}
	\tag{TR2}
	\begin{aligned}
		\text{Minimize} \quad & (g,u)_{L^2(0,T)} + \beta \TV(u)\\
		\text{such that} \quad & \norm{u-v}_{L^1(0,T)}\leq \Delta^k,\quad u \in  \Uad.
	\end{aligned}
\end{equation}
To obey the constraint $\norm{u-v}_{L^1(0,T)}\leq \Delta^k$, we introduce the 
so-called budget $B:=\lfloor \Delta^k/\Delta t \rfloor\in\N$. Notice that in the 
discrete scenario, we have
\[\norm{u-v}_{L^1(0,T)} = \Delta t \sum_{j=1}^n\abs{u^j-v^j}\leq \Delta^k,\] 
which means that $B\in\N$ marks an upper bound for the sum of all distances 
between each value at a time step of $u,v\in\Uad$.
Now, we define a value function which is slightly different from the previous 
one used for the proximal-gradient method for $b\in\{0,\dots,B\}$, $l\in\{1,\dots,d\}$ and $\iota\in\{1,\dots,n\}$ by setting
\begin{equation} 
	\label{eq:TR_value_func1}
	\Phi_{l, \iota, b}
	:=
	\min
	\set*{
		\Delta t\sum_{j=\iota}^{n} g^j \nu_{\kappa_j}+ \beta \sum_{j=\iota}^{n-1} \abs{\nu_{\kappa_{j+1}}-\nu_{\kappa_{j}}}
		\given
		\begin{aligned}
			& \kappa_\iota = l,\\
			& \kappa_{\iota+1},\dots,\kappa_n\in\{1,\dots,d\},\\
			& \sum_{j=\iota}^{n}\abs{\nu_{\kappa_j}-v^j} = b
		\end{aligned}
	},
\end{equation}
where we use the convention $\min \varnothing :=\infty$.
This means that if there do not exist $\kappa_{\iota+1},\dots,\kappa_n\in\{1,\dots,d\}$ with
\[b=\sum_{j=\iota}^{n}\abs{\nu_{\kappa_j}-v^j},\] 
 we have $\Phi_{l, \iota, b} :=\infty$.
 Imitating the arguments of the previous section, 
we see that
\[
	\forall l = 1, \ldots, d:
	\qquad
	\Phi_{l,n,\abs{\nu_l-v^n}} = \Delta t g^n \nu_l,
\]
while $\Phi_{l,n,b} = \infty$ for every pair $(l,b) \in \set{1,\ldots,d} \times \set{0,\dots,B}$ that cannot be 
represented as above.
In order to compute
$\Phi_{l,\iota,b}$ for $\iota < n$,
we have to minimize
\begin{equation*}
	\Delta t g^\iota \nu_{l}+ \beta \abs{\nu_{\kappa_{\iota+1}}-\nu_{l}}
	+
	\Delta t\sum_{j=\iota+1}^{n} g^j \nu_{\kappa_j}+ \beta \sum_{j=\iota+1}^{n-1} \abs{\nu_{\kappa_{j+1}}-\nu_{\kappa_{j}}}
\end{equation*}
w.r.t.\ $\kappa_{\iota+1}, \ldots, \kappa_n \in \set{1,\ldots,d}$
such that
$\sum_{j=\iota+1}^{n}\abs{\nu_{\kappa_j}-v^j} = b - \abs{\nu_l-v^\iota}$.

As in \cref{sect:prox_subproblems}, this allows to rewrite \eqref{eq:TR_value_func1} as 
\begin{equation}
	\label{eq:TR_value_func3}
	\Phi_{l, \iota, b}
	:=
	\min
	\set*{
		\Delta t g^\iota \nu_l
		+ \beta \abs{\nu_{\kappa_{\iota+1}}-\nu_l} + \Phi_{\kappa_{\iota+1},\iota+1,b - \tilde b}
		\given
		\begin{aligned}
			& \kappa_{\iota+1}\in\{1,\dots,d\},\\
			& \abs{\nu_l-v^\iota} = \tilde b \le b
		\end{aligned}
	}.
\end{equation}
Now, for every $l\in\{1,\dots,d\}$ and $b\in\{0,\dots,B\}$, we calculate 
$\Phi_{l,\iota,b}$ for $\iota = n-1,\dots,1$ while saving the corresponding 
minimizer in a structure $U\in\R^{d\times n-1\times B}$ given by
\begin{equation*}
	U_{l, \iota, b}
	:=
	\argmin
		\set*{
		\Delta t g^\iota \nu_l
		+ \beta \abs{\nu_{\kappa_{\iota+1}}-\nu_l} + \Phi_{\kappa_{\iota+1},\iota+1,b - \tilde b}
		\given
		\begin{aligned}
			& \kappa_{\iota+1}\in\{1,\dots,d\},\\
			& \abs{\nu_l-v^\iota} = \tilde b \le b
		\end{aligned}
	},
\end{equation*}
while using $\argmin\varnothing := 0$.
The pair $(l_1,b_1)$ minimizing $\Phi_{l,1,b}$
w.r.t.\ $l \in \set{1,\ldots,d}$, $b \in \set{0,\ldots,B}$
can be used to to reconstruct the solution $u$ by calculating the values $(l_\iota,b_\iota)$
for all $\iota\in\{1,\dots,n-1\}$ via
\[
	l_{\iota+1}
	=
	U_{l_{\iota},\iota,b_\iota},
	\quad
	b_{\iota+1} = b_{\iota} - \abs{\nu_{l_{\iota}} - v^{\iota}}
\]
and setting $u^\iota= \nu_{l_\iota}$ for every $\iota\in\{1,\dots,n\}$.
\begin{remark}
	\label{rem:trust-subproblem}
	Similarly to \cref{rem:prox-subproblem},
	only a $d\times 2\times B$ array $\Phi$ is needed when the above calculations are carried out.
	Here, for every time step of $\{1,\dots,n-1\}$,
	we have to test $d$ target values for $d$ possible settings of $l$ and at maximum $B$ possible values for $b$,
	suggesting that this procedure has a runtime of $\OO(d^2nB)$.
	Since $B$ is of order $n$
	(for trust-region radii $\Delta^k$ which are bounded from below and from above),
	this results in a total runtime of $\OO(d^2n^2)$.

	In contrast to the method for proximal-gradient subproblems,
	it is not possible to adapt the above procedure to
	general non-equidistant meshes since the definition of $B$ depends on the uniform mesh size $\Delta t$.
	However, in the important case that all occurring interval lengths $t_j - t_{j-1}$
	are integer multiples of a minimal length,
	it is possible to transfer the ideas.
\end{remark}

\section{Numerical examples}
\label{sec:numerics}
To study the properties and quality of the proximal-gradient (PG) and trust-region (TR) algorithm  using the Bellman principle, we consider a Lotka-Volterra fishing problem motivated by \cite[Chapter 4]{Sager2012} aswell as a signal reconstruction problem involving a convolution investigated in \cite{LeyfferManns2021}. 

The problems will be discretized using a grid with $n$ equidistant grid points, where we will test 
different values for $n$ ranging from $256$ to $4096$. For (PG), 
we will choose the algorithmic parameters $\eta = 10^{-6}$, $\theta = -\tfrac 12$ and 
$\tau^0 = 0.01$, while (TR) will be initiated with an initial trust-region
radius of $\Delta^0 = 0.4$ for the Lotka-Volterra problem and $\Delta^0 = 0.125$ 
for the signal reconstruction problem. The algorithms are implemented in Julia Version 1.6.3 and all results 
are calculated using an Intel(R) Core(TM) i9-10900 CPU @ 2.80GHz on a Linux OS.
\subsection{Lotka-Volterra fishing problem}
	\label{sec:LV}
	For parameters $\alpha_1$, $\alpha_2$, $\gamma_1$, 
	$\gamma_2$, $\theta_1$, $\theta_2$, $\beta$, $T>0$ and an initial state 
	$y_0\in\R^2$, the Lotka-Volterra fishing problem is given by
	\begin{equation}
		\label{eq:LV}
		\tag{LV}
		\begin{aligned}
			\text{Minimize} \quad & \frac 12\int_0^T (y_1(t) -1)^2 + (y_2(t)-1)^2 \, \d t + \beta\TV(u)\\
			\text{such that} \quad
			&y_1'(t) = \alpha_1 y_1(t) - \alpha_2 y_1(t) y_2(t) - \theta_1 y_1(t) u(t) && \text{a.e.\ on }(0,T)\\
			&y_2'(t) = \gamma_1 y_1(t) y_2(t) - \gamma_2 y_2(t) - \theta_2 y_2(t) u(t) && \text{a.e.\ on }(0,T)\\
			&y(0) = y_0,\quad u(t)\in\set{0,1}\text{ a.e.\ on } (0,T).
		\end{aligned}
	\end{equation}
As stated in \cite[Chapter 4]{Sager2012}, the problem does not admit a solution 
when the term $\TV(u)$ is not present. However, the optimal objective value 
can be approximated arbitrarily close when $u$ is switching often enough.

We can write \eqref{eq:LV} in the form of \eqref{eq:prob} by defining an operator 
$S\colon L^2(0,T)\to W^{1,1}(0,T,\R^2)$ mapping a function $u\in L^2(0,T)$ to the 
unique solution of the ordinary differential equation (ODE) in \eqref{eq:LV}. Thus, we have 
\[F(u) = \frac 12\int_0^T \paren*{S(u) - \begin{pmatrix}1\\1\end{pmatrix}}^\top \paren*{S(u) - \begin{pmatrix}1\\1\end{pmatrix}} \, \d t.\]
It can be verified that $F$ is bounded from below by 0 and continuous on 
$L^1(0,T)$ if $S$ is continuous. The continuity of $S$ together with its Fréchet differentiability
can be shown by employing the implicit function theorem,
see \cref{sec:Appendix}.
The derivative $S'(u)$ can be characterized with the adjoint equation corresponding to the ODE in 
\eqref{eq:LV}. In the implementation, we solved all occurring ODEs using the explicit
Euler method.

We tested the algorithms by using 1000 randomly generated initial guesses $u_0$
constructed such that $u_0$ switches values 32 times at uniformly chosen unique grid points from $\{1,\dots,n\}$, 
where the 33 corresponding control levels are picked randomly from $\{\nu_1,\dots,\nu_d\}$.
Also, we used the parameters
\[(\alpha_1,\alpha_2,\gamma_1,\gamma_2,\theta_1,\theta_2) = (1,1,1,1,0.4,0.2),\quad y_0 =\begin{pmatrix}
	0.5\\0.7
\end{pmatrix},\quad T = 12 \]
and $\beta=0.0001$
in \eqref{eq:LV}.

\begin{table}[htp]
	\centering
	\begin{tabular}{r|cc|cc|cc}
		$n$ & \multicolumn{2}{c|}{range of objectives} & \multicolumn{2}{c|}{average time [$s$]}     & \multicolumn{2}{c}{average iterations} \\
		& PG                  & TR                 & PG                  & TR                  & PG                & TR                 \\ \hline
		256  & $[0.738,4.126]$     & $[0.716,0.775]$    & $2.32\cdot 10^{-3}$ & $1.09\cdot 10^{-3}$ & 2.93              & 25.17              \\
		512  & $[0.749,3.522]$     & $[0.694,0.720]$    & $4.47\cdot 10^{-3}$  & $3.16\cdot 10^{-3}$ & 3.01              & 30.98              \\
		1024 & $[0.707,3.324]$     & $[0.683,0.704]$    & $8.39\cdot 10^{-3}$ & $8.83\cdot 10^{-3}$ & 2.79              & 34.78              \\
		2048 & $[0.708,3.188]$     & $[0.678,0.697]$    & $1.63\cdot 10^{-2}$ & $3.57\cdot 10^{-2}$ & 2.79              & 48.18              \\
		4096 & $[0.787,3.261]$     & $[0.675,0.694]$    & $3.31\cdot 10^{-2}$ & $3.80\cdot 10^{-1}$ & 2.64          		& 147.1             
	\end{tabular}
	\caption{Results of applying (PG) and (TR) 1000 times to \eqref{eq:LV} with random start point $u_0$ for different 
		grid sizes $n$.}
	\label{table_LV}
\end{table}
\begin{figure}[htp]
	\centering
	\includegraphics{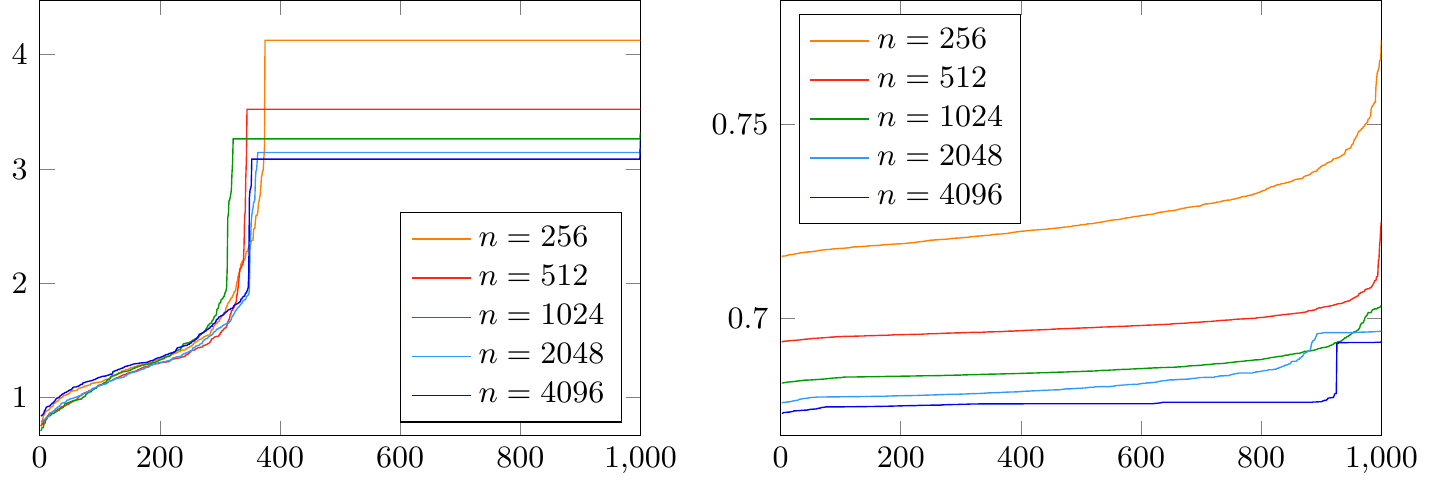}
	\caption{Distribution of 1000 objective values for \eqref{eq:LV} calculated by (PG) (left) and (TR) (right) for different choices of 
		$n$ and random start functions $u_0\in\Uad$.}
	\label{distribution_LV}
\end{figure}
\begin{figure}[htp]
	\includegraphics{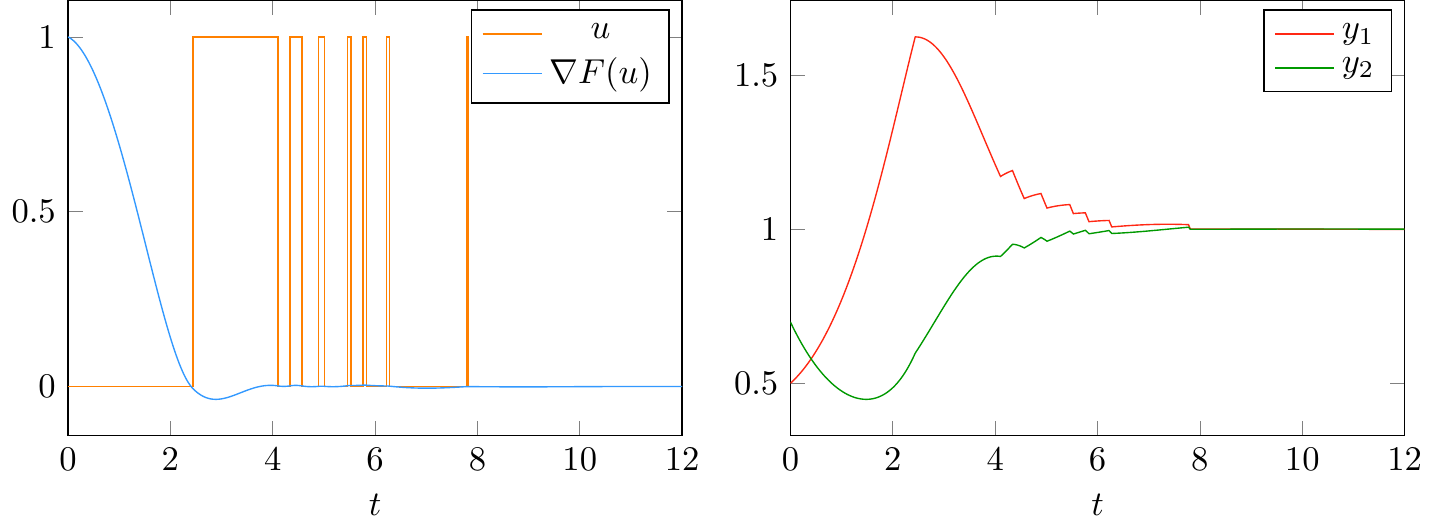}
	\caption{Solution of \eqref{eq:LV} gained with (TR) by iteratively enlarging grid until $n=4096$ with $\beta=10^{-4}$, objective: $0.6749$, time: $0.061s$.
	$\nabla F(u)$ is scaled such that $\norm{\nabla F(u)}_{L^\infty(0,T)}=1$.
}
	\label{plot_LV}
\end{figure}

In \cref{table_LV}, we can see that (TR) generally produces far better results than (PG) with 
comparable computing times. This may be 
due to the fact that (PG) is not suited for non-convex 
optimization problems. Indeed, in more than 50\% of all cases for every grid size, the 
solution generated by (PG) will be zero in every grid point after 2 iterations of the outer loop, 
which can be observed by 
interpreting the distributions of the objective values in \cref{distribution_LV} and the last column in \cref{table_LV}. 

The best results can be achieved by starting (TR) with a randomly generated start function $u_0$ 
on a grid of size $n=256$ and using the corresponding solution as a start function on a 
refined grid (with halved time step size), which will be repeated until arriving at $n=4096$. Indeed, using 
this method testing again 1000 randomly generated start functions, we arrive at an 
objective range of $[0.6749,0.6789]$ with an average computing time of 0.128s. 

Solutions as displayed in \cref{plot_LV} are 
competitive, since the optimal objective value for the relaxed problem (allowing 
$u(t)\in[0,1]$) without the total variation term (i.e., $\beta = 0$) is given by 0.67204,
cf. \cite[Chapter~4.1]{Sager2012}.
Note that $\nabla F(u)$ is equal or close to zero whenever $u$ switches.

\subsection{Signal reconstruction problem}
	\label{ex:convolution}
	To compare our results with the SLIP-method derived in 
	\cite{LeyfferManns2021}, we consider the problem
	\begin{equation}
		\label{eq:CV}
		\tag{SR}
		\begin{aligned}
		\text{Minimize} \quad & \frac 12 \norm{Ku - f}_{L^2(t_0,t_f)}^2 + \beta\TV(u)\\
		\text{such that} \quad &u(t)\in\set{-2,-1,0,1,2}\text{ a.e.\ on } (t_0, t_f),
		\end{aligned}
	\end{equation}
where $Ku:=k*u$ for the convolution kernel
\begin{equation*}
	k(t)
	:=
	-\frac{\sqrt{2}}{10}\chi_{[0,\infty)}(t)\omega_0
	\exp\paren*{-\frac{\omega_0 (t-1)}{\sqrt{2}}}
	\sin\paren*{\frac{\omega_0 (t-1)}{\sqrt{2}}}
.
\end{equation*}
Furthermore, we use the data $\omega_0= \pi$, $t_0=-1$, $t_f=1$ aswell as $f(t)
:= \tfrac 25\cos(2\pi t)$. In \cite[Proposition 5.1]{LeyfferManns2021}, it is
shown that $F(u):=\frac 12 \norm{Ku - f}_{L^2\paren{t_0,t_f}}^2$ is continuously
differentiable with $\nabla F(u) = K^*(Ku-f)$, where $K^*$ denotes the adjoint
operator of $K$.
Since the objective is bounded from below by zero, the problem meets our assumptions.
 
As described before, the problem will be discretized using a grid
$\{t_0,\dots,t_n\}$ with the equidistant mesh size $\Delta t:=\frac{t_f-t_0}{n}$
and setting $u(t):=\sum_{j=1}^n u^j\chi_{(t_{j-1},t_j)}(t)$.
We further introduce the vectors
$\mathpzc{u}=(u^1,\dots,u^n)^\top$,
$\mathpzc{f} = (f(t_0), f(t_1), \ldots, f(t_n))^\top$.
In this scenario,
the evaluation of the convolution $Ku$
in a grid point $t_i$, $i\in\{0,\dots,n\}$
can be calculated as a simple matrix-vector product:
Since
\begin{equation*}
	\begin{aligned}
(K u)(t_i) &= \int_{t_0}^{t_i}k(t_i -\tau)u(\tau)\,\d \tau = \sum_{j=1}^i	u^j\int_{t_{j-1}}^{t_j}k(t_i -\tau)\,\d \tau = \sum_{j=1}^i	u^j\int_{t_i-t_j}^{t_i-t_{j-1}}k(\tau)\,\d \tau,
	\end{aligned}
\end{equation*}
we can write $(K u)(t_i) = \parens*{\mathpzc{K}\mathpzc{u}}_{i+1}$ for 
$i=0,\dots,n$ with the matrix $\mathpzc{K} = (k_{lj})_{(l,j)\in I}$, 
$I=\set{1,\dots,n+1}\times\set{1,\dots,n}$ given by
\[k_{lj} = \begin{cases}
		\int_{t_{l-1}-t_{j}}^{t_{l-1}-t_{j-1}}k(\tau)\,\d \tau,&\text{if } j< l,\\
		0,&\text{if } j\geq l.
\end{cases}\]
Note that $\mathpzc{K}$ is a Toeplitz matrix with zeros on and above the main diagonal, thus it is only necessary to compute 
$\int_{t_{l-1}-t_1}^{t_{l-1}-t_0}k(\tau)\,\d \tau$
for $l=2,\dots,n+1$. This will be done using the 5th-order Gauß-Legendre quadrature rule.

In order to discretize the objective function,
we linearly interpolate the values $(K u)(t_i)$, i.e.,
we redefine
\[
(K u)(t) := \sum_{i=0}^{n} (K u)(t_i)\phi_i(t)= (\mathpzc{Ku})^\top\mathpzc{\phi}(t),
\quad
f (t):=\sum_{i=0}^{n} f(t_i)\phi_i(t)= \mathpzc f^\top\mathpzc{\phi}(t),
\]
where
$\phi_0, \phi_1, \ldots, \phi_n$ are the usual
(piecewise linear) hat functions on the 
grid $\set{t_0, t_1, \ldots, t_n}$
and $\phi(t) = (\phi_0(t),\dots,\phi_{n}(t))^\top$.
With this, the first part of the objective in \eqref{eq:CV} is discretized as
\[\frac 12 \int_{t_0}^{t_n}(\mathpzc{Ku-f})^\top \phi(t)\phi(t)^\top (\mathpzc{Ku-f})\,\d t = \frac 12 (\mathpzc{Ku-f})^\top M (\mathpzc{Ku-f}) \]
with $M:=(m_{ij})_{i,j=1}^{n+1}$, $m_{ij} = \int_{t_0}^{t_n}\phi_{i-1}(t)\phi_{j-1}(t)\,\d t$.
It is easy to see that the derivative
of the first part of the (discretized) objective
in this scenario is given by
\[\nabla F(u) = \mathpzc K^\top M (\mathpzc{Ku-f}).\]

\begin{figure}[p]
	\centering
	\includegraphics{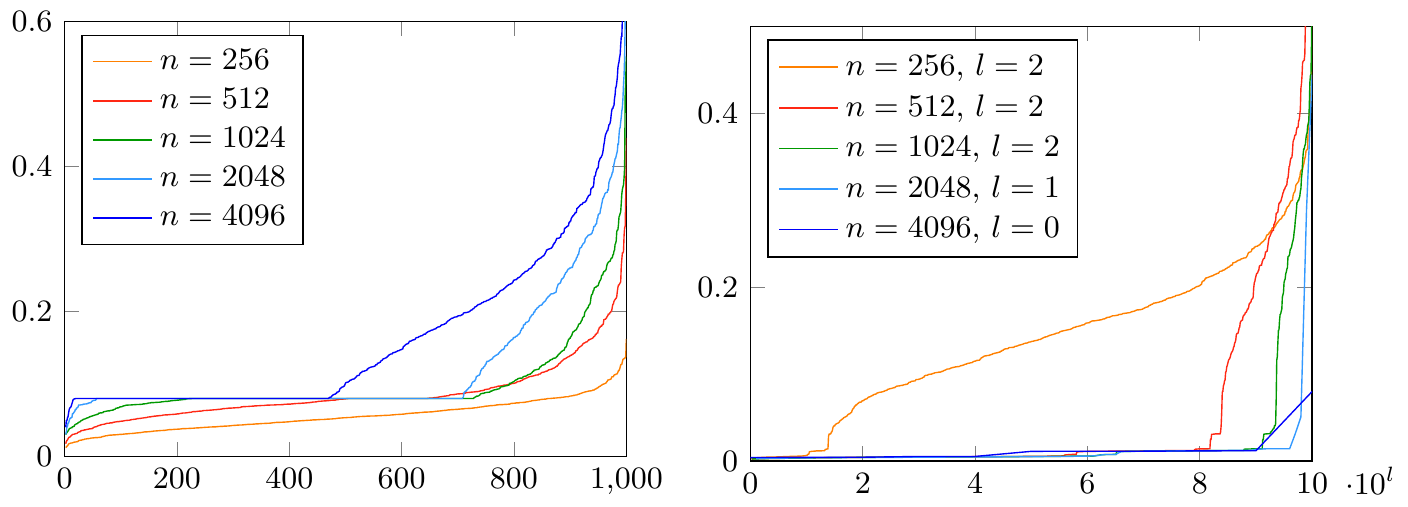}
	\caption{Distribution of objective values for \eqref{eq:CV} calculated by (PG) (left) and (TR) (right) for different choices of 
		$n$ and random start functions $u_0\in\Uad$.}
	\label{distribution_CV}
\end{figure}

Now, we tested different random start functions $u_0$ again, constructed as in \cref{sec:LV} but switching 128 times.
The results are displayed in \cref{table_CV} and \cref{distribution_CV}.
\begin{table}[p]
	\centering
	\begin{tabular}{r|cc|cc|cc}
		\multicolumn{1}{c|}{n}   & \multicolumn{2}{c|}{range of objectives} & \multicolumn{2}{c|}{average time [$s$]}     & \multicolumn{2}{c}{average iterations} \\
		& PG                  & TR                 & PG                  & TR                  & PG                & TR                 \\ \hline
		\multicolumn{1}{c|}{256} & $[0.0117,0.1620]$     & $[0.0032,0.5098]$    & $6.61\cdot 10^{-3}$ & $1.38\cdot 10^{-2}$  & 3.34   & 13.79   \\
		\multicolumn{1}{c|}{512} & $[0.0171,0.4396]$     & $[0.0025,0.6613]$    & $2.18\cdot 10^{-2}$  & $4.70\cdot 10^{-1}$ & 2.76   & 93.89 \\
		1024                     & $[0.0225,0.5518]$     & $[0.0024,0.5486]$    & $1.02\cdot 10^{-1}$ & $4.79$ & 2.42              & 126.94 \\
		2048                     & $[0.0276,0.6500]$     & $[0.0026,0.4441]$ & $7.29\cdot 10^{-1}$ & $57.1$ & 2.322              & 192.75     \\
		4096                     & $[0.0335,0.7340]$     & $[0.0040,0.0800]$ & $5.75 $ & $481.1$ & 2.43            & 236.09
	\end{tabular}
	\caption{Results of applying the (PG) and (TR) $10^l$ times to \eqref{eq:CV} with random start point $u_0$ for 
		different grid sizes $n$, where $l=2$ for every grid size when applying (PG), while 
		$l=2$ for $n\in\{256,512,1024\}$, $l=1$ for $n=2048$ and $l=0$ for $n=4096$ when 
		applying (TR).}
	\label{table_CV}
\end{table}
\begin{figure}[p]	
	\includegraphics{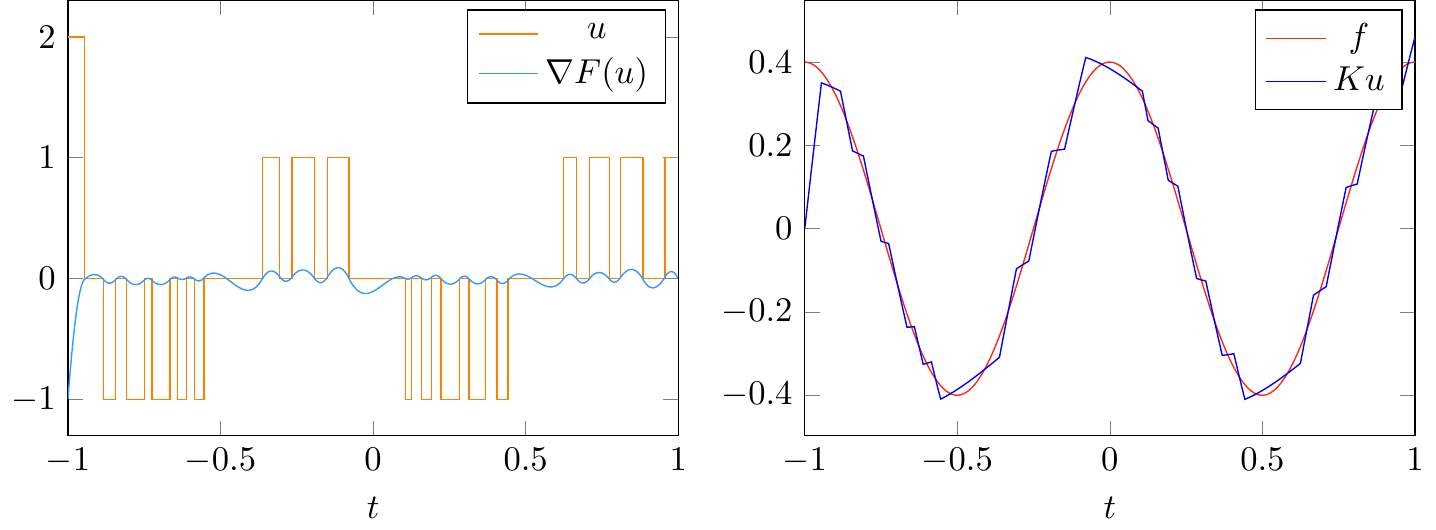}
	\caption{Solution of \eqref{eq:CV} gained with (TR) by iteratively enlarging grid until $n=4096$ with $\beta=10^{-4}$, objective: $2.04\cdot 10^{-3}$, time: $89.3s$.
	$\nabla F(u)$ is scaled such that $\norm{\nabla F(u)}_{L^\infty(0,T)}=1$.
}\label{plot_CV}
\end{figure}

Once more, (PG) performs worse than (TR), where better results are achieved for small grid sizes $n$. Again, in a lot of cases, the solution will be zero at every grid point, as the distributions tend to stagnate in a certain objective value in \cref{distribution_CV}.

On the other hand, (TR) behaves as expected, with larger grid sizes resulting 
in (generally) smaller objective values with a higher average computing time and iteration 
number. This motivates to again 
refine the grid starting with a random start point $u_0$ on the grid size $n=256$ 
until $n=4096$, such that when testing 50 random start functions, we arrive at an 
objective range of $[2.04\cdot 10^{-3},1.35\cdot 10^{-2}]$ and an average computing time 
of $261.4s$. A good solution is showcased in \cref{plot_CV}. Note that in 
\cite[Chapter 5]{LeyfferManns2021}, the presented solution was calculated on a grid of 
size $n=2048$ with an objective value of $4.339\cdot 10^{-3}$ in $1.698\cdot 10^4s$.

Varying the number of jumps for start functions generated as in \cref{sec:LV} 
will have a noticeable impact on the quality of the received solutions, even
when refining the grid. 
The algorithms were also tested using other randomizations for $u_0$. 
For example, when assigning a random value of $\{\nu_1,\dots,\nu_d\}$ to $u_0$ in every grid point,
the results gained by (PG) will be worse (compared to \cref{table_LV}, \cref{distribution_LV}) for \eqref{eq:LV} and in a lot of cases zero in every grid point for \eqref{eq:CV}. 
On the other hand, (TR) is able to generate comparable solutions 
in this scenario,
where the results get slightly better for \eqref{eq:CV} and slightly worse for \eqref{eq:LV}.

\subsection{Runtime of trust-region subproblem solver}
In \cite{SeverittManns2022}, two methods to solve trust-region subproblems discretized as 
a shortest path problem on a directed acyclic graph were tested. One method used a 
topological sorting of the nodes (TOP), while the other arises from the Dijkstra 
algorithm using a heuristic which gives a lower bound for the cost to reach the sink from 
any node in the graph (Astar). In order to compare these methods to our solver derived 
with the Bellman principle (BP), we will test it using instances of \eqref{eq:TR_subprob2} where
\[
\{\nu_1,\dots,\nu_d\} = \{-2,\dots,23\}
,\quad
T = 1
,\quad
\Delta t = \tfrac1n
,\quad
\Delta^k=\tfrac18
,\quad
B = \tfrac n8
\]
and $\beta\in[0,1]$, $v^j\in\{\nu_1,\dots,\nu_d\}$ are chosen uniformly, while $g^j$ is 
chosen from a normal distribution with mean 0 and variance 1 for all $j=1,\dots,n$. These 
instances are constructed in a way to resemble the problem (SH) from \cite[Section 
5.1]{SeverittManns2022}.
Note that the runtime of (BP)
does not depend on the values of $g$ and $v$,
since the main work is to evaluate \eqref{eq:TR_value_func3}
and its effort is independent of $g$ and $v$.

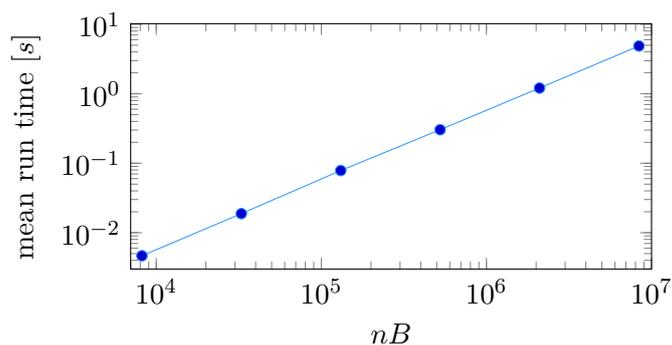
\begin{figure}[htp]	
	\centering
	\begin{tikzpicture}[>=latex]
		\begin{loglogaxis}[xlabel={$nB$}, ylabel={mean run time [$s$]},
			y = .4cm,
			xmin = 7000, xmax = 1e7, ymin = .003
			]
			\addplot+[nablacol] file {data_files/runtimes.dat};
		\end{loglogaxis}
	\end{tikzpicture}
	\caption{Mean run times of the trust-region subproblem solver (BP) over the product of the grid size $n$ and the budget $B$ for randomly generated instances of \eqref{eq:TR_subprob2}.}
	\label{fig:runtimes}
\end{figure}

For every choice of $n\in\{2^8,\dots,2^{13}\}$, we will employ (BP) 20 times, each time with a different randomization. The mean run times will be displayed in correspondence to the value $nB = \tfrac{n^2}{8}$ in \cref{fig:runtimes} to compare our results with those of \cite[Figure 3]{SeverittManns2022}.

We can see that the runtime of (BP) depends linearly on $nB$, which is not surprising since the expected runtime scales linearly with respect to this product,
see \cref{rem:trust-subproblem}. For small values of $nB$, i.e. close to $10^4$, (BP) seems to be faster 
than (Astar) and slightly slower than (TOP), while for large $nB$, i.e. close to $10^7$, it appears that (BP) has 
roughly the same runtime as (Astar). However, since our computational setup is different than the one used in 
\cite{SeverittManns2022}, these comparisons should be taken with a grain of salt.

Since proximal-gradient subproblems can be solved faster than trust-region subproblems, we tried to develop a mixed algorithm, where a trust-region step instead of a proximal-gradient step will be done whenever $u_{k+1} = u_k$. However, this did not yield satisfactory results.

\section{Conclusion and outlook}
\label{sec:concl}
We investigated first and second order optimality conditions for integer control optimization problems using a switching point reformulation. The essential tool to show these conditions was the full representation of a piecewise constant function, allowing only switches between adjacent control levels. Non-local optimality conditions involving back-and-forth switches were also derived. 

Next, we showed convergence results of a proximal-gradient algorithm and used the Bellman principle to efficiently solve the corresponding subproblems. This method was adapted for subproblems of a trust-region method suggested in \cite{LeyfferManns2021}.

Testing the algorithms on two numerical examples showed that the proximal-gradient algorithm is not able to produce satisfactory results, while the trust-region method will give a good solution in most cases. Given that the best solutions found for our problems still do not meet the necessary optimality conditions derived in \cref{subsec:local_opt_con}, it may be advantageous to optimize the location of the switching points of such a solution with second-order methods by using the derivatives of \cref{thm:second_derivatives};
combined with the insertion and removal of switches
by utilizing \cref{thm:non_optimal_switch,thm:non_optimal_non_switch}.

Furthermore, the runtime of the subproblem solver could be improved by adapting the ideas from \cite{SeverittManns2022}. To be more precise, when given a heuristic to estimate a lower bound for the cost of a path in $U$, it may be possible to reduce the number of calculations carried out.

In a lot of applications, multiple decisions interact with a system simultaneously, motivating a generalization of the ideas presented in this paper for multidimensional control functions. 

\appendix
\section{Solution operator of the Lotka-Volterra ODE}
\label{sec:Appendix}
We prove that the operator $S$ introduced in \cref{sec:LV} is well defined and Fréchet differentiable.
To this end, we define $e\colon W^{1,1}(0,T,\R^2)\times L^1(0,T)\to L^1(0,T,\R^2)\times \R^2$ via
\[e(y,u) = \begin{pmatrix}
	y' - f(y,u) \\
	y(0) - y_0
\end{pmatrix}\quad\text{with}\quad f(y,u) = \begin{pmatrix}
	\alpha_1 y_1 - \alpha_2 y_1 y_2 - \theta_1 y_1 u\\
	\gamma_1 y_1 y_2 - \gamma_2 y_2 - \theta_2 y_2 u
\end{pmatrix}\]
and employ the implicit function theorem.
In order to show the Fréchet differentiability of $e$,
we only have to verify that
$f \colon W^{1,1}(0,T,\R^2) \times L^1(0,T) \to L^1(0,T,\R^2)$
is Fréchet differentiable, since all the other terms are linear and bounded.
First, we expect that the partial derivatives of $f$
are given by
\begin{equation*}
	f_y(y,u)z = 
	\begin{pmatrix}
		\alpha_1 z_1 - \alpha_2(y_1z_2+z_1y_2) \\
		\gamma_1(y_1z_2 + z_1y_2)-\gamma_2 z_2
	\end{pmatrix}
	,\quad
	f_u(y,u)v
	=
	\begin{pmatrix}
		- \theta_1 y_1v \\
		- \theta_2 y_2v
	\end{pmatrix}
	.
\end{equation*}
Now, the remainder is given by
\begin{equation*}
	f(y+z,u+v) - f(y,u)
	-
	f_y(y,u)z + f_u(y,u)v
	=
	\begin{pmatrix}
		-\alpha_2 z_1z_2 - \theta_1z_1v \\
		\gamma_1z_1z_2 - \theta_2z_2v
	\end{pmatrix}.
\end{equation*}
Now, it is true that
\begin{align*}
	&\frac{\norm*{
			\begin{pmatrix}
				-\alpha_2 z_1z_2 - \theta_1z_1v\\
				\gamma_1z_1z_2 - \theta_2z_2v
			\end{pmatrix}
	}_{L^1(0,T,\R^2)}}
	{\norm*{
			\begin{pmatrix} z\\v\end{pmatrix}
	}_{W^{1,1}(0,T,\R^2)\times L^1(0,T)}} 
	=
	\frac{
		\norm{\alpha_2 z_1z_2 + \theta_1z_1v}_{L^1(0,T)}
		+
		\norm{\gamma_1z_1z_2 - \theta_2z_2v}_{L^1(0,T)}
	}{
		\norm{z_1}_{W^{1,1}(0,T)}
		+
		\norm{z_2}_{W^{1,1}(0,T)}
		+
		\norm{v}_{L^1(0,T)}
	}
	\\
	&\leq C \frac{\norm{z_1}_{L^1(0,T)}\norm{z_2}_{L^\infty(0,T)} + \norm{z_1}_{L^\infty(0,T)}\norm{v}_{L^1(0,T)}+\norm{z_2}_{L^\infty(0,T)}\norm{v}_{L^1(0,T)}}
{
		\norm{z_1}_{W^{1,1}(0,T)}
		+
		\norm{z_2}_{W^{1,1}(0,T)}
		+
		\norm{v}_{L^1(0,T)}
	}
	\\
	&\leq
	C_2
	\frac
{
		\norm{z_1}_{W^{1,1}(0,T)}^2
		+
		\norm{z_2}_{W^{1,1}(0,T)}^2
		+
		\norm{v}_{L^1(0,T)}^2
	}
{
		\norm{z_1}_{W^{1,1}(0,T)}
		+
		\norm{z_2}_{W^{1,1}(0,T)}
		+
		\norm{v}_{L^1(0,T)}
	}
	\to 0
\end{align*}
if $(z, v) \to 0$
in
$W^{1,1}(0,T,\R^2) \times L^1(0,T)$.
Here, we used Hölder's inequality and
the continuous embeddings $W^{1,1}(0,T)\hookrightarrow L^\infty(0,T)$, $L^\infty(0,T)\hookrightarrow L^1(0,T)$.
Thus, $f$ and $e$ are Fréchet differentiable.
Moreover, the partial derivative of $e$ w.r.t.\ $y$
is given by
\[e_y(y,u)z = \begin{pmatrix}
	z' - f_y(y,u)z\\z(0)
\end{pmatrix} \in L^1(0,T,\R^2)\times \R^2,\]
which is
a linear ordinary differential operator with an initial condition.
Thus, the continuous invertibility of $e_y(y,u)$
follows from
\cite[Chapter~5, Theorem~1.3]{GajewskiGroegerZacharias1975}.
Using the implicit function theorem, we see that $S$ is well defined aswell as Fréchet 
differentiable. 

\printbibliography

\end{document}